\documentclass[12pt]{article}
\usepackage{amsmath,latexsym,amssymb}
\usepackage{eucal}
\newtheorem{thm}{Theorem}[section]
\newtheorem{lem}[thm]{Lemma}
\newtheorem{prop}[thm]{Proposition}
\newtheorem{df}[thm]{Definition}
\newtheorem{exam}[thm]{Example}

\newcommand{\id}{\mathrm{id}}

\newcommand{\cA}{\mathcal{A}}
\newcommand{\cB}{\mathcal{B}}
\newcommand{\cC}{\mathcal{C}}
\newcommand{\cD}{\mathcal{D}}
\newcommand{\cE}{\mathcal{E}}
\newcommand{\cM}{\mathcal{M}}

\newcommand{\cCf}{\mathcal{C}}
\newcommand{\cCinf}{\mathcal{C}}

\newcommand{\hlambda}{\Lambda}

\newcommand{\hateta}{{\eta_G}}

\newcommand{\fH}{\mathfrak{H}}

\DeclareMathOperator{\Ad}{Ad}

\DeclareMathOperator{\Aut}{Aut}

\DeclareMathOperator{\End}{End}

\DeclareMathOperator{\Hom}{Hom}
\DeclareMathOperator{\Irr}{Irr}

\begin{document}

\title{
Relative center construction for $G$-graded C$^*$-tensor categories and Longo-Rehren inclusions
}
\author{Toshihiko MASUDA\footnote{The author is supported by JSPS 22K03341} \\
Faculty of Mathematics, Kyushu University \\
744, Motooka, Nishi-ku, Fukuoka 819-0395, Japan \\
E-mail: masuda@math.kyushu-u.ac.jp}

\date{}

\maketitle

\begin{center}
\textit{
Dedicated to Professor Masaki Izumi on the occasion of his 60th birthday} 
\end{center}

\begin{abstract}
Gelaki-Naidu-Nikshych and Turaev-Virelizier showed the existence of $G$-braiding on the
        relative Drinfeld center of a $G$-graded tensor category. We will explain this concept
        from the viewpoint of Longo-Rehren inclusions. 
\end{abstract}

\section{Introduction}\label{sec:intro}
In the theory of tensor categories, braided tensor categories have been extensively studied,
because they play important roles in various situations, e.g., conformal field theory, 
topological quantum field theory, and so on.
An important method to obtain braided tensor categories is 
the Drinfeld center construction. 
An operator algebraic counterpart of the center construction is
asymptotic inclusions introduced by A. Ocneanu.
He called it  quantum double construction,
and clarified properties of the quantum double construction in terms of TQFT
by introducing a notion of Tube algebras for tensor categories. 
In particular, he showed the existence of braiding.
See \cite{EK-book} for the detail of his theory.

There are other construction which corresponds to the Drinfeld center
in operator algebras. One construction is Popa's symmetric enveloping algebra \cite{Po-symm},
and another construction is Longo-Rehren inclusion \cite{LR}. From the viewpoint of tensor categories
Longo-Rehren inclusion is suitable for the study of operator algebraic aspect of center construction.
(See \cite{M-asym}, \cite{M-LR} on the relation of these constructions. )

In \cite{Iz-LR} and \cite{Iz-LR-II}, M. Izumi studied the structure of Longo-Rehren inclusion.
He introduced the notion of  
half-braiding, and proved Ocneanu's statement for asymptotic inclusion without use of TQFT, and
present many interesting examples. Since then, the notion of a half braiding has been widely used
in the study of the Drinfeld center.

The Drinfeld center construction can be generalized to relative version. Namely,
for given inclusion of C$^*$-tensor categories, we have a natural notion of 
the relative Drinfeld center for this inclusion. (See \cite{Kw-Dcenter} for instance.)

As a special case, 
we can define a relative center for a $G$-graded tensor category for a discrete group $G$.
In \cite{GNN-center}, Gelaki-Naidu-Nikshych investigate the structure of 
$G$-graded fusion category for a finite group $G$. They showed the existence of $G$-braiding 
on the relative Drinfeld center. This result is generalized by Turaev-Virelizier 
\cite{TuVi-center} for tensor categories with pivotal structure and a discrete group $G$. 
The construction in \cite{GNN-center} is conceptual and
one in \cite{TuVi-center} is rather constructive. 

In this paper, we explain this result in the framework of operator algebras. 
Namely, we concretely realize the relative center construction and  
$G$-braidings in the set of  endomorphisms of a Longo-Rehren 
inclusion. To construct $G$-braiding, we need to construct an action of $G$ on the relative center.
We first show that the existence of an action by constructing an action of $G$ on a Longo-Rehren inclusion.
Then we construct a $G$-braiding. 
Our construction of $G$-braiding is similar to the construction of braiding in AQFT. 
We also discuss the relation between relative centers and relative Tube algebras,
which is a generalization of Ocneanu's statement.

This paper is organized as follows. In \S \ref{sec:LoRe}, we recall the definition of Longo-Rehren inclusion. 
Our presentation is based on \cite{M-LR} rather than the original definition given in \cite{LR},
because the definition in \cite{LR} is valid only for C$^*$-fusion categories.
When a given C$^*$-tensor category is $G$-graded, we show the existence of an action of  $G$, which
reflects the $G$-grading structure. 

In \S \ref{sec:gradingcenter}, we study the relative Drinfeld center and the endomorphism categories 
of Longo-Rehren inclusions by using a notion of a half braiding as in \cite{Iz-LR}. 
The main purpose is the construction of a $G$-braiding on the relative center.
Our construction of a $G$-braiding is similar to the construction of braidings in AQFT.

In \S \ref{sec:full}, we study the relationship between the relative Drinfeld  center 
for a $G$-graded C$^*$-tensor category and  the  
Drinfeld center for a full category as in \cite[\S 3B]{GNN-center}.
We show that a Drinfeld center for a full category is isomorphic
to the $G$-equivariant tensor category of a relative center. Key of our argument is 
that the Longo-Rehren inclusion for a full category is obtained by the crossed product construction.

In \S \ref{sec:Gcenter}, we apply the results in \S \ref{sec:LoRe} for C$^*$-tensor category with
group actions. In this case, the existence of a $G$-braiding is proved rather easily than
that given in \S \ref{sec:gradingcenter}.

In Appendix, we present detail of M{\"u}ger's remark \cite[Remark 5.1]{Mu-subII} on
Tube algebras.

\section{Preliminaries on Longo-Rehren inclusion}\label{sec:LoRe}
Our standard references are \cite{Tak-book} for operator algebras, 
\cite{BKLR-Qsys-book}, \cite{EK-book}, \cite{Nes-Tu} 
for subfactor theory and tensor categories, and 
\cite{Iz-fusion},  \cite{Kosaki-seoul} 
for sector theory.

\subsection{Endomorphisms of type III factors}\label{subsec:endo}

Let $\cM$ be a factor of type III, 
$(L^2(\cM), L^2(\cM)_+, J)$ the standard form of $\cM$,
$\End(\cM)$ the set of all endomorphisms, and 
$\End_0(\cM)$ the set of all endomorphisms with finite index.

Two endomorphisms $\pi,\rho\in \End(\cM)$ are said to be equivalent if $\pi=\Ad u\circ \rho$ for some
unitary $u\in \cM$, and we write $\pi \sim \rho$.

Let $\{\rho_i\}_{i\in I}\subset \End(\cM)$ be an at most countable set, and take 
a family of isometries $\{w_i\}_{i\in I}\subset \cM$ with $\sum_{i\in I}w_iw_i^*=1$.
A direct sum of $\{\rho_i\}_{i\in I}$ is defined by  
$\sum_{i\in I}w_i\rho_i(x)w_i^*\in \End(\cM)$, whose equivalence class does not 
depend on the choice of $\{w_i\}_{i\in I}$. 
When $I$ is an infinite set, this sum converges in the strong topology.

For $\rho,\sigma\in \End(\cM)$, their intertwiner space is defined by $(\rho,\sigma):=
\{T\in \cM\mid T\rho(x)=\sigma(x)T \mbox{ for all } x\in \cM\}$. When $\rho$ is irreducible,
i.e., $(\rho,\rho)=\rho(\cM)'\cap \cM=\mathbb{C}1$, $(\rho,\sigma)$ has an inner product by
$\langle T,S\rangle 1=S^*T$.

For $\rho\in \End(\cM)$,  ${}_\rho L^2(\cM)$ (resp. $L^2(\cM)_\rho$) denotes an $\cM$-$\cM$ bimodule 
defined by $a\cdot \xi \cdot b:=\rho(a)Jb^*J \xi$ (resp. 
$a\cdot \xi \cdot b:=aJ\rho(b)^*J \xi$). 
We define a conjugate endomorphism $\bar{\rho}$ of $\rho$ by  an $\cM$-$\cM$ bimodule isomorphism
${}_\rho L^2(\cM)\cong L^2(\cM)_{\bar{\rho}}$.
For $\sigma\in \End_0(\cM)$, 
we call a pair $(R_\sigma,\bar{R}_\sigma)$ a solution of a conjugate equation for $(\sigma,\bar{\sigma})$ if 
\[
 R_\sigma\in (\id,\bar{\sigma}\sigma),\,
 \bar{R}_\sigma\in (\id,\sigma\bar{\sigma}), \,
R_\sigma^*R_\sigma=\bar{R}_\sigma^*\bar{R}_\sigma=d(\sigma), \,R_\sigma^* \bar{\sigma}(\bar{R}_\sigma)=
\bar{R}_\sigma^* {\sigma}({R}_\sigma)=1.
\]

Fix a faithful normal state $\psi_0$  on $\cM$.  
Let $\rho \in \End_0(\cM)$, and $(R_\rho, \bar{R}_\rho)$ be a solution of a conjugate equation for $\rho$.
Let $\phi_\rho(x):=R_{\rho}^*\bar{\rho}(x)R_\rho$ be a (non-normalized) standard left inverse.
Define a standard implementation $u_\rho\in B(L^2(\cM))$ of $\rho$ by 
$u_\rho (a\psi_0^{\frac{1}{2}}):=\rho(a)(\psi_0\circ \phi_\rho)^{\frac{1}{2}}$.
Here $\psi_0^{\frac{1}{2}}\in L^2(\cM)_+$ is a vector implementing $\psi_0$.
We also use the notation $u_\rho(v)=\rho(v)$ for $v\in L^2(\cM)$.
We have $Tu_\rho=JT^*Ju_{\rho'}$, $T\in (\rho,\rho')$. In other words, 
$T\rho(v)=\rho'(v)T$, $v\in L^2(\cM)$.
It is easy to see that  $u_\rho^*$ is given by $R_\rho^* JR_\rho^* Ju_{\bar{\rho}}$, i.e. $u_\rho^*(v)=R_\rho^* \bar{\rho}(v)R_\rho$.

\subsection{Longo-Rehren inclusion}\label{subsec:LR}

Let $\cC\subset \End_0(\cM)$ be a rigid C$^*$-tensor category, and
$\Irr(\cC)$  the set of all representatives of irreducible objects of $\cC$.

\smallskip

\noindent
\textbf{Notation.} In this note, we often take summation over $\mathrm{Irr}(\cC)$.
To simplify notation, we denote this summation by $\sum_{\pi\in \cC}$.
(We believe no confusion arises.)

\smallskip

Let $j:\cM\rightarrow \cM^{\mathrm{op}}$ the canonical conjugate linear isomorphism
given by $j(x)=x^*$, and set $\cA:=\cM\otimes \cM^{\mathrm{op}}$.
For $\pi\in \End(\cM)$, let $\pi^{\mathrm{op}}:=j\circ \pi\circ j^{-1}\in \End(\cM^{\mathrm{op}})$.

Fix an orthonormal basis $\{T_{\pi,\rho}^{\sigma,e}\}_{e=1}^{N_{\pi,\rho}^\sigma} \subset (\sigma,\pi\rho)$,
$\pi,\rho,\sigma\in \mathrm{Irr}(\cC)$, $N_{\pi,\rho}^\sigma=\dim(\sigma,\pi\rho)$.
Let
\[
\tilde{T}_{\pi,\rho}^\sigma:=\sum_{e=1}^{N_{\pi,\rho}^\sigma}
T_{\pi,\rho}^{\sigma,e}\otimes
j(T_{\pi,\rho}^{\sigma,e}) \in (\sigma\otimes \sigma^{\mathrm{op}},\pi\rho\otimes \pi^{\mathrm{op}}\rho^{\mathrm{op}})
\]
Due to the conjugate linearity of $j$, $\tilde{T}_{\pi,\rho}^\sigma$ is well-defined  
and the following lemma can be verified. 
\begin{lem}\label{lem:cano-inter}
We have the following. 
\begin{align*}
\sum_{\sigma\in {\cC}}
\tilde{T}_{\xi,\eta}^\sigma
\tilde{T}_{\sigma,\zeta}^\pi&=\sum_{\sigma \in {\cC}} \xi\otimes  \xi^{\mathrm{op}}
(\tilde{T}_{\eta,\zeta}^\sigma)\tilde{T}_{\xi,\sigma}^\pi,\,\,\,  
\tilde{T}_{\xi,\eta}^{\zeta*}\tilde{T}_{\xi,\eta}^{\zeta'}=\delta_{\zeta,\zeta'}
N_{\xi,\eta}^\zeta, \\
T_{\xi,\eta}^{\zeta*}&=\frac{d(\zeta)}{d(\xi)}\zeta\otimes \zeta^\mathrm{op}(R_{\eta\otimes \eta^{\mathrm{op}}}) ^*
T_{\zeta, \bar{\eta}}^\xi 
=\frac{d(\zeta)}{d(\eta)}\bar{R}_{\xi\otimes \xi^\mathrm{op}}^* \xi \otimes \xi^{\mathrm{op}}
(T_{\bar{\xi},\zeta}^\eta).
\end{align*}
\end{lem}

Let $\cA^{0}(\cC)$ a $*$-algebra generated by $\cA$ and $\{\lambda_\pi\}_{\pi\in \Irr(\cC)}$ by the following relation;
\[
 \lambda_\pi a=\pi\otimes \pi^{\mathrm{op}}(a)\lambda_\pi,\,\, a\in \cA , \,\,
\lambda_\pi\lambda_\rho=\sum_{\sigma\in \cC}\tilde{T}_{\pi,\rho}^\sigma \lambda_\sigma,\,\,
\lambda_\pi^*=R^*_{\pi\otimes \pi^{\mathrm{op}}}\lambda_{\bar{\pi}}.
\]
Here $R_{\pi\otimes \pi^{\mathrm{op}}}=R_\pi\otimes j(R_\pi)$.
Thus any element of $\cA^0(\cC)$ is expressed by $\sum_{\pi\in \Irr(\cC)}a_\pi \lambda_\pi$, where 
$a_\pi \in \cA$, and $a_\pi=0$ except for finitely many $\pi$.

It is convenient to extend $\lambda_\sigma$ for any $\sigma\in \cC$ as follows.
Let $\sigma(x)=\sum_{k}w_{k}\pi_k(x)w_{k}^*$, $\pi_k\in \Irr(\cC)$,
 be an irreducible decomposition, and define
$\lambda_\sigma:=\sum_{k}w_{k}\otimes j(w_{k})\lambda_{\pi_k}$. 
Then $\lambda_\sigma$ is well-defined, i.e., it does not depend on the choice of $\{w_k\}$.
The above relation is rewritten as follows;
\[
\lambda_\pi a=\pi\otimes \pi^{\mathrm{op}}(a)\lambda_\pi, \,\,
 \lambda_\pi\lambda_\rho=\lambda_{\pi\rho}, \,\, 
\lambda_\pi^*=R_{\pi\otimes \pi^\mathrm{op}}^*\lambda_{\bar{\pi}}, \,\, 
(T\otimes 1)\lambda_\pi=(1\otimes j(T^*))\lambda_\rho 
\]
for $\pi,\rho\in \cC$, $T\in (\pi,\rho)$.
 In particular, we have $\lambda_\rho^*\lambda_\rho=(R_\rho^*\otimes j(R_\rho^*))\lambda_{\bar{\rho}\rho}=
1\otimes j(R_\rho^*R_\rho)=d(\rho)$.

Let  $\fH:=\bigoplus_{\cC}L^2(\cA)=\ell^2(\Irr(\cC), L^2(\cA))$. 
Define 
$ \gamma_L(a), \gamma_R(a), \kappa_\pi, \gamma_L(\lambda_\pi), \gamma_R(\lambda_\pi)
\in B(\fH)$, $a\in \cA$, 
$\pi\in \Irr(\cC)$, and a conjugate linear operator $\tilde{J}:\fH\rightarrow \fH$ by
\[
 (\gamma_L(a)v)(\pi):=(\pi\otimes \pi^{\mathrm{op}})(a)v(\pi),\,\,\,  (\gamma_R(a)v)(\pi):=Ja^*Jv(\pi)=v(\pi)a,  
\] 
\[
(\gamma_L(\lambda_\pi) v)(\rho):=\sum_{\sigma\in \cC}\sqrt{\frac{d(\rho)}{d(\sigma)}}\tilde{T}_{\rho,\pi}^{\sigma} v(\sigma),
( \kappa_\pi v)(\rho):=\sum_{\sigma\in \cC}\sqrt{\frac{d(\rho)}{d(\sigma)}}\tilde{T}_{\pi,\rho}^{\sigma}v(\sigma),
\gamma_R(\lambda_\pi):=u_{\pi \otimes \pi^{\mathrm{op}}}^*\kappa_\pi,
\]
\[
  (\tilde{J}v)(\pi):=d(\pi)J\tilde{T}_{\pi,\bar{\pi}}^{\id*}{\pi}(v(\bar{\pi}))=
J\bar{R}_{\pi\otimes \pi^{\mathrm{op}}}^*{\pi}(v(\bar{\pi})).
\]
Here $u_{\pi\otimes \pi^{\mathrm{op}}}$ is the standard implementation of $\pi\otimes \pi^{\mathrm{op}}$.

\begin{lem}\label{lem:LR-relation}
$(1)$ We have following relations;
\begin{align*}
(\pi\otimes \pi^{\mathrm{op}})(\kappa_\rho)\kappa_\pi&=   \sum_{\sigma\in \cC}\tilde{T}_{\pi,\rho}^\sigma\kappa_\sigma, \,\,
\kappa_\pi^*=\bar{R}^*_{\pi\otimes  \pi^\mathrm{op}} \pi\otimes \pi^{\mathrm{op}} (\kappa_{\bar{\pi}}),  \\
\kappa_\pi\gamma_L(a\lambda_\rho)&=\pi\otimes \pi^{\mathrm{op}}(\gamma_L(a\lambda_\rho))\kappa_\pi.
\end{align*}
$(2)$ By $\gamma_L$ and $\gamma_R$, we can define an $\cA^0(\cC)$-bimodule structure on $\fH$, that is, 
$\gamma_L$ (resp. $\gamma_R$) defines a left (resp. right) action of $\cA^0(\cC)$
with $\gamma_L(\cA^0(\cC))\subset \gamma_R(\cA^0(\cC))'$.  \\
$(3)$ $\tilde{J}$ is an anti-unitary involution and $\tilde{J}\gamma_L(a^*)\tilde{J}=\gamma_R(a)$,\,\, $a\in \cA^0(\cC)$.
\end{lem}
\textbf{Proof.} The statement (1) and the fact $\gamma_L$ being a left action follow from Lemma \ref{lem:cano-inter}.
Then the rest of the statement (2)  follows from the statement (1) and Lemma \ref{lem:cano-inter}. 
The statement (3) follows from  Lemma \ref{lem:cano-inter} and the property of canonical implementation of endomorphisms.
\hfill $\Box$

\medskip

\begin{df}\label{df:LRinclusion}
Let  $\cB:=\gamma_L(\cA^0(\cC))''$.
An inclusion $\gamma_L(\cA)\subset \cB$ is called 
the Longo-Rehren inclusion associated with $\cC$.
\end{df}

Like the crossed product construction, there exists a conditional expectation $E:\cB\rightarrow \gamma_L(\cA)$
by $E(a\lambda_\pi)=a\delta_{\id,\pi}$. 
For $a\in \cB$, let $a_\pi:=E(a\lambda_\pi^*)$.
Then we can expand $a=\sum_{\pi\in \cC}d(\pi)a_\pi \lambda_\pi$ uniquely.
(Let  $\psi:=\psi_0\otimes \psi^{\mathrm{op}}_0\circ E$.
This expansion converges with respect to the norm $\|\cdot \|_{\psi}$.)
In the following, we omit $\gamma_L$ and 
denote the  Longo-Rehren inclusion by $\cA\subset \cB$.
Namely, the Longo-Rehren inclusion $\cA\subset \cB$ is characterized by 
(1) $\cB$ is a  von Neumann algebra generated by $\cA$ and $\{\lambda_\pi\}_{\pi\in \Irr(\cC)}$ satisfying
the above relation, (2) 
there exists  a conditional expectation $E:\cB\rightarrow \cA$ given 
by $E(a\lambda_\pi)=a\delta_{\id,\pi}$.
In what follows, we regard $\gamma_L$ (resp. $\gamma_R$) as a left (resp. right) representation of $\cB$ on $\fH$.
(Note Lemma \ref{lem:LR-relation} (3) implies $\gamma_R(\cB)=\tilde{J}\gamma_L(\cB)\tilde{J}$.)

We can easily verify $\cA'\cap \cB=\mathbb{C}1$ by using expansion of elements in $\cB$. 
In particular $\cB$ is a factor.
Define $\delta_0\in \fH$ by $\delta_0(\pi)=\delta_{\id,\pi}\psi^{\frac{1}{2}}\otimes (\psi^{\mathrm{op}})^{\frac{1}{2}}$. Then
we can identify $(\fH,\gamma_L, \delta_0)$ with a GNS  representation of $\cB$ by $\psi$.

\subsection{Longo-Rehren inclusion for a $G$-graded category}\label{subsec:LRgrade}

Let $G$ be a discrete group. The neutral element of $G$ is denoted by $0$, although 
the operation of $G$ is written multiplicatively.
From now on, we assume $\cC$ is $G$-graded, and $\cC_g$  
denotes the $g$-component for $g\in G$. 
We can decompose $\fH$ as 
\[
\fH=\bigoplus_{g\in G}\fH_g,\,\,\,
\fH_g:=\bigoplus_{\pi\in \cC_g}L^2(\cA)=\ell^2(\Irr(\cCf_g), L^2(\cA)).
\]

Let $\cB_g$ be  $\sigma$-weak closure of 
$\{\sum_{\pi\in \mathrm{Irr}(\cC_g)}a_\pi \lambda_\pi\mid a_\pi\in \cA\}\subset \cB$.
Note that $\cB_0$ is a subfactor of $\cB$, and it is a Longo-Rehren inclusion for $\cC_0$.
We have $\gamma_L(\cB_g)\fH_h\subset \fH_{h{g^{-1}}}$, $\gamma_R(\cB_g)\fH_h\subset \fH_{g^{-1}h}$.
In particular $\fH_g$ is a $\cB_0$-$\cB_0$-bimodule.

In what follows, we denote by $\Hom(_{\cB_0} X_{\cB_0}, {}_{\cB_0} Y_{\cB_0})$
the set of all $\cB_0$ bilinear map for $\cB_0$-$\cB_0$ bimodules ${}_{\cB_0} X_{\cB_0}$,
${}_{\cB_0} Y_{\cB_0}$, and
$\End({}_{\cB_0} X_{\cB_0})=
\Hom({}_{\cB_0} X_{\cB_0}, {}_{\cB_0} X_{\cB_0})$.

Set
\[
 \fH^{\mathrm{bdd}}_{g,\cB_0}:=\{v\in \fH_g\mid \mbox{there exists }C>0 \mbox{ such that }
\|\gamma_R(a)v\|\leq C\psi(aa^*) \mbox{ for all }a\in \cB_0\}.
\]
For $v \in \fH^{\mathrm{bdd}}_{g,\cB_0}$, 
a map $L_\psi(v)$ given by 
$L_\psi(v)\gamma_R(a)\delta_0=\gamma_R(a)v$ 
extends to one in  $B(\fH_0, \fH_g)$. 
Any  element of $\fH^{\mathrm{bdd}}_{g,\cB_0}$ is of the form $\gamma_L(a)\delta_0$ for some $a\in \cB_{g^{-1}}$, 
\cite[Lemma IV.3.3]{Tak-book}. (Note $\delta_0\in \fH_0$.)
Therefore 
\[
L_\psi(\gamma_L(a)\delta_0)\gamma_R(b)\delta_0=\gamma_R(b)\gamma_L(a)\delta_0
=\gamma_L(a)\gamma_R(b)\delta_0
\]
holds, and hence we get $L_\psi(\gamma_L(a)\delta_0)=\gamma_L(a)$.
\begin{prop}\label{B0-auto}
 There exists an action $\beta_g$ of $G$ on $\cB_0$ such that 
$\fH_g\cong {}_{\beta_g} \fH_0$ as  $\cB_0$-$\cB_0$-bimodules.
\end{prop}
\textbf{Proof.} 
By the definition of an inner product on a relative tensor product \cite[Proposition 3.5]{Tak-book}
and $L_\psi(\gamma_L(a)\delta_0)=\gamma_L(a)$, we have
\begin{align*}
 \langle \gamma_L(a)\delta_0\otimes_\psi v, \gamma_L(b)\delta_0\otimes_\psi w\rangle&=\langle 
L_\psi(\gamma_L(b)\delta_0)^*L_\psi(\gamma_L(a)\delta_0)v,w\rangle 
=\langle \gamma_L(a)v,\gamma_L(b)w\rangle .
\end{align*}

Thus a map 
$W_{g,h}\in \Hom(\fH_g\otimes_\psi
\fH_h, \fH_{hg}) $
given by 
\[ 
W_{g,h}(\gamma_L(a)\delta_0\otimes_\psi w)
=\gamma_L(a)w 
\]
is an isometry.

Let $w\in \fH_{hg}$. Fix $\pi_0\in \cC_{g}$, and  
put $v:=d(\pi_0)^{-1}\gamma_L(\lambda_{\pi_0})w\in \fH_{h}$.
Then 
\[
 W_{g,h}(\gamma_L(\lambda_{\pi_0}^*)\delta_0\otimes_\psi v)=d(\pi_0)^{-1}\gamma_L(\lambda_{\pi_0}^*\lambda_{\pi_0}) w=w.
\]
Thus  $W_{g,h}$ is surjective, and hence 
$W_{g,h}$ is a unitary. 
Since $\cB_0$ is a type III factor, there exists $\beta_k\in \End(\cB_0)$ 
such that ${}_{\cB_0}\fH_{k,\cB_0} \cong {}_{\cB_0,\beta_k}\fH_{0,\cB_0}$ 
as $\cB_0$-$\cB_0$-bimodules.
Take a unitary  $U_k 
\in \Hom({}_{\cB_0}\fH_{k,\cB_0},{}_{\cB_0,\beta_k}\fH_{0,\cB_0})
$. 
Since $U_k$ commutes with $\gamma_R(\cB_0)$, 
$U_k$ is in $ \gamma_L(\cB_k)$. 
Let 
$I \in \Hom( {}_{\beta_h}\fH_0\otimes_\psi{}_{\beta_g}\fH_0, {}_{\beta_g\beta_h}\fH_0 )$
be a natural unitary given by
$\gamma_L(a)\delta_0\otimes_\psi \gamma_L(b)\delta_0\rightarrow \gamma_L(\beta_g(a)b)\delta_0$.
Then   $ U_{gh}W_{h,g}(U_h^*\otimes_\psi U_g^*) I^*:
{}_{\beta_g\beta_h}\fH_0\rightarrow  {}_{\beta_{gh}}\fH_0
$ is a $\cB_0$-$\cB_0$ bimodule isomorphism.
This  implies $\beta_g\beta_h\sim \beta_{gh}$ and  
$\beta_g\in \Aut(\cB_0)$. 
Here 
\begin{align*} 
U_{gh}W_{h,g}(U_h^*\otimes_\psi U_g^*)I^* \gamma_L(a)\delta_0&=
U_{gh}W_{h,g}(U_h^*\otimes_\psi U_g^*)(\delta_0\otimes_\psi  \gamma_L(a)\delta_0) \\ &=
U_{gh}W_{h,g}(U_h^* \delta_0\otimes_\psi U_g^*\gamma_L(a)\delta_0)  =
U_{gh}U_h^* U_g^*\gamma_L(a)\delta_0.
\end{align*}
Hence  we obtain $U_{gh}W_{h,g}(U_h^*\otimes_\psi U_g^*)I^*=U_{gh}U_h^* U_g^*$.
In particular, $U_{gh}U_h^* U_g^*\in \gamma_R(\cB_0)'=\gamma_L(\cB_0)$.
Let $u(g,h)\in \cB_0$ be a unitary by 
$\gamma_L(u(g,h))=U_{gh}U_h^*U_g^*$.
It is easy to see that $(\beta_g, u(g,h))$ is a cocycle twisted action. Hence  we can assume 
$U_gU_h=U_{gh}$ by the second cohomology vanishing theorem \cite{Su-homoI}, and $\beta_g$ is an action on $\cB_0$.
\hfill$\Box$

\medskip

\noindent
\textbf{Remark.} When $|\Irr(\cC)|<\infty$, we can show the existence of an action $\beta$ by 
 applying results in  \cite[\S 3]{BKLR-Qsys-book} without use of bimodules.

\section{Relative center construction for a C$^*$-tensor category with grading}\label{sec:gradingcenter}
From now on, we discuss the relative center $Z_{\cC_0}(\cC)$.
We will collect facts on the relative center based on \cite{Iz-LR}, \cite{Kw-Dcenter}.

In this paper, we treat countable infinite direct sum of elements of $\cC$ when $\Irr(\cC)$ is an infinite set,
i.e., ind-object in the sense of \cite{Nesh-Yama-DrinRep}.
Abusing notation, we write $\sigma \in \cC$ even in such a case.

At first, we recall the notion of half braiding. (See \cite[Definition 4.2]{Iz-LR}, \cite[Definition 2.1]{Kw-Dcenter}.)

\begin{df}\label{df:halfbraiding}
Let $\sigma\in \cCinf$. A family of unitaries $\{\mathcal{E}_\sigma(\pi)\}_{\pi\in \cC_0}$ is
a half braiding for $\sigma$  if  \\
$(1)$ $\mathcal{E}_\sigma(\pi)\in (\sigma\pi,\pi\sigma)$, \\
$(2)$ $\cE_\sigma(\pi\rho)=\pi(\cE_\sigma(\rho))\cE_\sigma(\pi)$,\,\, $\pi,\rho\in \cC_0$,  \\
$(3)$ $\cE_\sigma(\rho)\sigma(T)=T \cE_\sigma(\pi)$,\,\, $T\in (\pi,\rho)$.
\end{df}

Before stating the following lemma, we remark on the notation.
Let  $a\in \gamma_R(\cA)'=\cA\otimes B(L^2(\cA)) (\subset B(\fH_0))$. 
Then we can consider $(\sigma\otimes \id)\otimes \id_{B(L^2(\cA))}(a)$ for $\sigma\in \End(\cM)$.
In the following we often denote this by $\sigma\otimes \id(a)$, or  
$\sigma(a)$ for simplicity.

\begin{lem}\label{lem:halfglobal}
Let $\{\cE_\sigma(\pi)\}_{\pi\in \cC_0}$ be a half braiding for $\sigma$.
Define  a unitary $\cE_\sigma\in \End(\fH_{0,\cA})$  
by $(\cE_\sigma v)(\xi)=(\cE_\sigma(\xi)\otimes 1)v(\xi)$.
 Then the following statements hold. \\
$(1)$ $\cE_\sigma\sigma(\gamma_L(a))=\gamma_L(\sigma\otimes \id(a))\cE_\sigma$, $a\in \cA$. \\
$(2)$ $\gamma_L(\lambda_\pi)\cE_\sigma =\gamma_L(\cE_\sigma(\pi)\otimes 1)\cE_\sigma\sigma\gamma_L(\lambda_\pi)$, 
$\pi\in \Irr(\cC_0)$.  \\
$(3)$ $\kappa_\pi\cE_\sigma =\pi(\cE_\sigma)(\cE_\sigma(\pi)\otimes 1)\sigma(\kappa_\pi)$,
$\pi\in \Irr(\cC_0)$. \\
(See the above remark about the meaning of $\sigma\gamma_L(\cdot)$.)

Conversely, if a unitary 
$\cE_\sigma\in \End(\fH_{0,\cA})$  
satisfies (1) and (2), or (1) and (3), 
then $\cE_\sigma$ comes from a half braiding as above.
\end{lem}
\textbf{Proof.}
(1) For $a\in \cA$,
\begin{align*}
\left( \cE_\sigma\sigma\gamma_L(a)v \right)(\xi)&=
(\cE_\sigma(\xi)\otimes 1)(\sigma\xi\otimes \xi^{\mathrm{op}})(a)v(\xi) \\
&=
(\xi\sigma\otimes \xi^{\mathrm{op}})(a)(\cE_\sigma(\xi)\otimes 1)v(\xi) 
=\left(\gamma_L(\sigma\otimes \id(a))\cE_\sigma v\right)(\xi).
\end{align*}

(2) By the half braiding property 
 \[
       T \cE_\sigma(\eta)=\xi(\cE_\sigma(\pi))\cE_\sigma(\xi)\sigma(T),\,\,\, T\in (\eta,\xi\pi),
      \]
we have
\begin{align*}
& \left(\gamma_L(\lambda_\pi)\cE_\sigma v\right)(\xi) 
= \sum_{\eta\in \cC_0}\sqrt{\frac{d(\xi)}{d(\eta)}}T_{\xi,\pi}^\eta \left(\cE_\sigma v\right)(\eta) 
= \sum_{\eta\in \cC_0}\sqrt{\frac{d(\xi)}{d(\eta)}}
T_{\xi,\pi}^\eta (\cE_\sigma(\eta)\otimes 1)v(\eta) \\
&=
\sum_{\eta\in \cC_0} \sqrt{\frac{d(\xi)}{d(\eta)}}
\left(\xi(\cE_\sigma(\pi))\cE_\sigma(\xi)\otimes 1\right)\sigma(T_{\xi,\pi}^\eta)v(\eta) \\
&=
( \xi(\cE_\sigma(\pi))\otimes 1)
(\cE_\sigma\sigma \gamma_L(\lambda_\pi) v)(\xi) 
=
\big(
\gamma_L(\cE_\sigma(\pi)\otimes 1)\cE_\sigma\sigma\gamma_L(\lambda_\pi)
 v\big)(\eta).
\end{align*}

(3) In a similar way, we have
\begin{align*}
& \left(\kappa_\pi\cE_\sigma v\right)(\xi) 
= \sum_{\eta\in \cC_0}\sqrt{\frac{d(\xi)}{d(\eta)}}T_{\pi,\xi}^\eta \left(\cE_\sigma v\right)(\eta) 
= \sum_{\eta\in \cC_0}\sqrt{\frac{d(\xi)}{d(\eta)}}
T_{\pi,\xi}^\eta (\cE_\sigma(\eta)\otimes 1)v(\eta) \\
&=
\sum_{\eta\in \cC_0} \sqrt{\frac{d(\xi)}{d(\eta)}}
\left(\pi(\cE_\sigma(\xi))\cE_\sigma(\pi)\otimes 1\right)\sigma(T_{\pi,\xi}^\eta)v(\eta) \\
&=
( \pi(\cE_\sigma(\xi))\otimes 1)(\cE_\sigma(\pi)\otimes 1)
(\sigma (\kappa_\pi) v)(\xi) 
=
\big( \pi(\cE_\sigma)(\cE_\sigma(\pi)\otimes 1)\sigma(\kappa_\pi)
 v\big)(\xi) .
\end{align*}

Suppose a unitary 
$\cE_\sigma\in \End(\fH_{0,\cA})$  
satisfies (1). 
Then we can easily see $(\cE_\sigma v)(\xi)=(\cE_\sigma(\xi)\otimes 1)v(\xi)$ for some unitary 
$\cE_\sigma(\xi)\in (\sigma\xi,\xi\sigma)$. If we examine the proof of (2) and (3) above, we can see
the condition (2) or (3) is equivalent to the half braiding property.
\hfill$\Box$ 

\medskip

In the following,  we call the pair $(\sigma,\mathcal{E}_\sigma)$ a half braiding.

\begin{df}\label{df:eta-ext}
$(1)$ For a half braiding $(\sigma,\mathcal{E}_\sigma)$,
we define an $\eta$-extension $\eta(\sigma,\mathcal{E}_\sigma)\in \End(\cB_0)$ by
\begin{align*}
\eta(\sigma,\mathcal{E}_\sigma)(a)&=(\sigma\otimes \id)(a), \,\, a\in \cA, \\
\eta(\sigma,\mathcal{E}_\sigma)(\lambda_\pi)&=\left(\mathcal{E}_\sigma(\pi)^*\otimes 1\right)
\lambda_{\pi}, \,\, \pi\in \mathrm{Irr}(\cC_0). 
\end{align*}
$(2)$ 
Define the relative Drinfeld center   $Z_{\cC_0}(\cC)$ by 
\[
Z_{\cC_0}(\cC):= \{\theta\in \End(\cB_0)\mid \theta\sim \eta(\sigma,\mathcal{E}_\sigma) \mbox{ for some half braiding } (\sigma,\mathcal{E}_\sigma) \}. 
\]
\end{df}

Indeed, Lemma \ref{lem:halfglobal} says that 
$\cE_\sigma \sigma\gamma_L(a)\cE_\sigma^*=\gamma_L\left(\eta(\sigma,\cE_\sigma)(a)\right)$ for $a\in \cB_0$, and 
Definition \ref{df:eta-ext} (1) is justified.
We also remark that $\eta(\sigma,\cE_\sigma)(\lambda_\pi)=\left(\mathcal{E}_\sigma(\pi)^*\otimes 1\right)\lambda_\pi$
holds for any $\pi\in \cC_0$.

In a similar way, we can also define $\eta^{\mathrm{op}}(\sigma, \mathcal{E}_\sigma)$
as an extension of $\id\otimes \sigma^{\mathrm{op}}$. Namely,
\[
\eta^\mathrm{op}(\sigma,\mathcal{E}_\sigma)(\lambda_\pi)=(1\otimes j(\mathcal{E}_\sigma(\pi)^*))\lambda_\pi.
\]
By definition, $\eta(\sigma,\mathcal{E}_\sigma)$ and $\eta^{\mathrm{op}}(\rho,\mathcal{E}_\rho)$ commute,
and this fact will be important later. (See comment after Proposition \ref{prop:Vunitary}, and Theorem \ref{thm:Gbraid1}.)

As in \cite{Iz-LR}, \cite{Kw-Dcenter},
we have
\begin{align*}
Z_{\cC_0}(\cC)&=
\{\theta\in \End_0(\cB_0)\mid \theta\iota\sim \iota(\sigma\otimes \id) \mbox{ for some }\sigma\in \cC \} \\
&=\{\theta \in \End_0(\cB_0)\mid \theta \prec \iota(\sigma\otimes \id)\bar{\iota} \mbox{ for some }\sigma\in \cC\}, 
\,\,  (\mbox{when }|\mathrm{Irr}(\cC)|<\infty ).
\end{align*}

To simplify notations, we often denote
$\eta(\sigma, \mathcal{E}_\sigma)$, $\eta^{\mathrm{op}}(\sigma, \mathcal{E}_\sigma)$
by $\tilde{\sigma}$, $\tilde{\sigma}^{\mathrm{op}}$, respectively.
(Note that $\sigma$ can have several half braidings. So this notation may cause
some confusion. However, in this article, we think that such confusion will not occur.)

We recall properties of $\eta$-extension.
The following proposition is proved in \cite[Theorem 4.6]{Iz-LR}, \cite[Theorem 2.4]{BEK-LR}.
(Also see \cite{Lu.Ha-Dcenter-pre}.)
\begin{prop}\label{prop:eta-inter} 
Let $(\sigma, \mathcal{E}_\sigma)$, $(\rho,\mathcal{E}_\rho)$ be
half braidings. \\
$(1)$
\[
\left(\eta(\sigma,\mathcal{E}_\sigma), \eta(\rho,\mathcal{E}_\rho)
\right)=\left\{T\otimes 1\mid T\in (\sigma,\rho), \mathcal{E}_\rho(\pi)T=\pi(T)\mathcal{E}_\sigma(\pi)
\mbox{ for all }\pi\in \cC_0\right\}.
\]
$(2)$ Let $(R_\sigma,\bar{R}_\sigma)$ be a solution of conjugate equation for $\sigma$.
Define
\[
\mathcal{E}_{\bar{\sigma}}(\pi):=R_\sigma^*\bar{\sigma}(\mathcal{E}_\sigma(\pi)^*)\bar{\sigma}\pi(\bar{R}_\sigma).
\]
Then $(\bar{\sigma},\mathcal{E}_{\bar{\sigma}})$ be a half braiding, and
$\overline{\eta(\sigma,\mathcal{E}_\sigma)}\sim \eta(\bar{\sigma},\mathcal{E}_{\bar{\sigma}})$.
(See the proof about the meaning of the above formula for $d(\sigma)=\infty$.)
\end{prop}
\textbf{Proof.}  
We only present a proof of (2), which is slightly different from one in
\cite[Theorem 4.6(4)]{Iz-LR}.  

Before starting proof,  we explain the meaning of (2) for $\sigma$ with $d(\sigma)=\infty$.
For a moment, we assume $d(\sigma)<\infty$.
Let $\sigma(x)=\sum_{k}w_k\xi_k(x)w_k^*$,
$\bar{\sigma}(x)=\sum_{k}\bar{w}_k\bar{\xi}_k(x)\bar{w}_k^*$ be  irreducible decomposition.
Define 
$\cE_\sigma(\pi)_{i, j}$ and $\cE_{\bar{\sigma}}({\pi})_{i,j}$ by
\[
\cE_\sigma(\pi)_{i, j}:=\pi(w_j^*)\cE_\sigma(\pi)w_{i}\in (\xi_i \pi,\pi \xi_j),
\,\,\,
\cE_{\bar{\sigma}}({\pi})_{i,j}=
R_{\xi_i}^*\bar{\xi}_i(\cE_{\sigma}({\pi})_{i,j}^*)\bar{\xi}_i{\pi}(\bar{R}_{\xi_j})
\in (\bar{\xi}_i {\pi},{\pi} \bar{\xi}_j).
\]
Then we have
\[
\cE_\sigma(\pi)=\sum_{i,j}\pi(w_j)\cE_\sigma(\pi)_{i,j}w_i^*,\,\,
R_\sigma^*\bar{\sigma}(\mathcal{E}_\sigma(\pi)^*)\bar{\sigma}\pi(\bar{R}_\sigma)
=\sum_{i,j}{\pi}(\bar{w_j})\cE_{\bar{\sigma}}({\pi})_{i,j}\bar{w}_i^*.
\]
Since the right hand side of the second equality
makes sense for $\sigma$ with $d(\sigma)=\infty$,
one should understand the formula in (2) as this equality.

At first we assume $d(\sigma)<\infty$.
We begin with the following observation.
Let $u\in (\xi,\zeta)$, $\xi,\zeta \in \End(\cM)$, be a unitary. (Hence $\xi$ and $\zeta$ are equivalent.)
Then
$R_\xi^*\bar{\xi}(u^*)\bar{\xi}(\bar{R}_\zeta)=\bar{\zeta}(\bar{R}_\xi)^*\bar{\zeta}(u^*)R_\zeta$
and this is a unitary in $(\bar{\xi},\bar{\zeta})$.
(See \cite[Theorem 2.2.21]{Nes-Tu}.)

Choose $R_{\bar{\pi}\sigma}$ and $\bar{R}_{\sigma\bar{\pi}}$ as  
$R_{\bar{\pi}\sigma}=\bar{\sigma}(\bar{R}_{{\pi}})R_\sigma$, 
$\bar{R}_{\sigma\bar{\pi}}={\sigma}({R}_{{\pi}})\bar{R}_\sigma$.
By the
above observation,
\[
\bar{\sigma}{\pi}(\bar{R}_{\sigma\bar{\pi}})^*\bar{\sigma}{\pi}(\mathcal{E}_\sigma(\bar{\pi})^*)R_{\bar{\pi}\sigma}
\in ({\pi}\bar{\sigma},\bar{\sigma}{\pi})
\]
is a unitary. 
Together with  the half braiding condition
\[
\bar{R}_{{\pi}}={\pi}(\mathcal{E}_\sigma(\bar{\pi}))\mathcal{E}_\sigma({\pi})\sigma(\bar{R}_{{\pi}}),
\]
we have
\begin{align*}
\bar{\sigma}{\pi}(\bar{R}_{\sigma\bar{\pi}})^*\bar{\sigma}{\pi}(\mathcal{E}_\sigma(\bar{\pi})^*)R_{\bar{\pi}\sigma}
&=
\bar{\sigma}{\pi}(\bar{R}_\sigma^*) \bar{\sigma}{\pi}\sigma({R}_{{\pi}}^*)
\bar{\sigma}{\pi}(\mathcal{E}_\sigma(\bar{\pi})^*)
\bar{\sigma}(\bar{R}_{{\pi}})R_\sigma \\
&=
\bar{\sigma}{\pi}(\bar{R}_\sigma^*) \bar{\sigma}{\pi}\sigma({R}_{\pi}^*)
\bar{\sigma}\left(\mathcal{E}_\sigma({\pi})\sigma(\bar{R}_{\pi})\right)R_\sigma \\
&= \bar{\sigma}{\pi}(\bar{R}_\sigma^*) \bar{\sigma}(\mathcal{E}_{\sigma}({\pi}))
\bar{\sigma}\sigma{\pi}({R}_{{\pi}}^*)
\bar{\sigma}\sigma(\bar{R}_{{\pi}})R_\sigma \\
&= \bar{\sigma}{\pi}(\bar{R}_\sigma^*) \bar{\sigma}(\mathcal{E}_{\sigma}({\pi}))R_\sigma.
\end{align*}
This shows that $\mathcal{E}_{\bar{\sigma}}(\pi)$ defined in (2) is a unitary.
It is easy to see that
$\mathcal{E}_{\bar{\sigma}}(\pi)$ is a half braiding for $\bar{\sigma}$.

When $d(\sigma)=\infty$, we can show the same result by the formal computation. 
To present a rigorous proof, one should take the following approach.
Let $\cE_\sigma(\pi)_{i,j}$, and $\cE_{\bar{\sigma}}({\pi})_{i,j}$ 
be as above.
By 
\[
\mathcal{E}_\sigma({\pi})^*
=\bar{R}_{{\pi}}^*\pi(\mathcal{E}_\sigma(\bar{\pi}))\pi\sigma({R}_{{\pi}})
,\,\,
 R_{\bar{\pi}\xi}=\bar{\xi}(\bar{R}_{{\pi}}) R_\xi,\,\, 
 \bar{R}_{\xi\bar{\pi}}={\xi}(\bar{R}_{\pi}) \bar{R}_\xi,
\]
we can see 
\[
\cE_{\bar{\sigma}}({\pi})_{i,j}=
R_{\xi_i}^*\bar{\xi}_i(\cE_{\sigma}({\pi})_{i,j}^*)\bar{\xi}_i{\pi}(\bar{R}_{\xi_j})=
R_{\bar{\pi}\xi_i}^*\bar{\xi}_i\pi(\cE_\sigma(\bar{\pi})_{j,i}\bar{R}_{\xi_j\bar{\pi}}).
\]

As the above observation, 
we can verify  that
\[
 \sum_{i,j}{\pi}(\bar{w}_j)\cE_{\bar{\sigma}}({\pi})_{i,j}\bar{w}_i^*=
 \sum_{i,j}{\pi}(\bar{w}_j)
R_{\bar{\pi}\xi_i}^*\bar{\xi}_i\pi(\cE_\sigma(\bar{\pi})_{j,i}\bar{R}_{\xi_j\bar{\pi}})
\bar{w}_i^*
\in (\bar{\sigma}\pi, {\pi}\bar{\sigma})
\] 
is a unitary. It is easy to verify this is a half braiding for $\bar{\sigma}$.

We will show
$\overline{\eta(\sigma,\mathcal{E}_\sigma)}\sim \eta(\bar{\sigma},\mathcal{E}_{\bar{\sigma}})$.
When $d(\sigma)<\infty$, 
$(R_\sigma\otimes 1, \bar{R}_\sigma \otimes 1)$ is a solution of conjugate equation
for $\eta(\sigma,\mathcal{E}_\sigma)$ by (1).
Hence $\overline{\eta(\sigma,\cE)}=\eta(\bar{\sigma},\cE_{\bar{\sigma}})$.
However, this proof does not work in the case $d(\sigma)=\infty$.
Hence we will show the equivalence
${}_{\tilde{\sigma}} L^2(\cB_0) \cong L^2(\cB_0)_{\tilde{\bar{\sigma}}}$ of $\cB_0$-$\cB_0$-bimodules 
directly.
To motivate our argument, we will present a proof for $\sigma$ with $d(\sigma)<\infty$ at first.
We can easily see  that  $R_\sigma^*u_{\bar{\sigma}}$ is a unitary which
gives an equivalence
of ${}_\sigma L^2(\cM)$ and $L^2(\cM)_{\bar{\sigma}}$ as $\cM$-$\cM$-bimodules.
We will show a unitary $(R_\sigma^*u_{\bar{\sigma}}\otimes 1)\cE_\sigma^*$ 
gives a desired equivalence. 
For a left action of $\cB_0$, we have 
\begin{align*}
 (R_\sigma^*u_{\bar{\sigma}}\otimes 1)\cE_\sigma^*\gamma_L(\tilde{\sigma}(a))&=
 (R_\sigma^*u_{\bar{\sigma}}\otimes 1)\sigma\gamma_L(a)\cE_\sigma^* =
\gamma_L(a) (R_\sigma^*u_{\bar{\sigma}}\otimes \id)\cE_\sigma^* , \,\, a\in \cB_0.
\end{align*}

We next study a right action of $\cB_0$.
Since $\cE_\sigma$ commutes with $\gamma_R(\cA)$, 
\[
 (R_\sigma^*u_{\bar{\sigma}}\otimes 1)\cE_\sigma^*
\gamma_R(b)=
 (R_\sigma^*u_{\bar{\sigma}}\otimes 1)\gamma_R(b)
\cE_\sigma^*=
\gamma_R(\bar{\sigma}\otimes \id(b))  (R_\sigma^*u_{\bar{\sigma}}\otimes 1)\cE_\sigma^* 
\]
for $b\in \cA$.

By Lemma \ref{lem:halfglobal} (3),
\begin{align*}
& (R_\sigma^*u_{\bar{\sigma}}\otimes 1)\cE_\sigma^*\gamma_R(\lambda^{*}_\pi) 
= (R_\sigma^*u_{\bar{\sigma}}\otimes 1)\cE_\sigma^*\kappa_\pi^*u_{\pi\otimes \pi^{\mathrm{op}}} \\
&=  (R_\sigma^*u_{\bar{\sigma}}\otimes 1)
\sigma(\kappa_\pi^*)(\cE_\sigma(\pi)^*\otimes 1)\pi(\cE_\sigma^*)u_{\pi\otimes \pi^{\mathrm{op}}} 
=  \kappa_\pi^*(R_\sigma^*u_{\bar{\sigma}}\otimes 1)(\cE_\sigma(\pi)^*\otimes 1)u_{\pi\otimes \pi^{\mathrm{op}}}\cE_\sigma^* \\
&=  \kappa_\pi^*(R_\sigma^*\otimes 1) \bar{\sigma}(\cE_\sigma(\pi)^*\otimes 1)(u_{\bar{\sigma}}\otimes 1)
u_{\pi\otimes \pi^{\mathrm{op}}}\cE_\sigma^* \\
& =  \kappa_\pi^*(\pi(R_\sigma^*)\otimes 1) (\cE_{\bar{\sigma}}(\pi)\otimes 1)(u_{\bar{\sigma}}\otimes 1)u_{\pi\otimes \pi^{\mathrm{op}}}\cE_\sigma^* \\
&=  \kappa_\pi^*(\pi(R_\sigma^*)\otimes 1) (J\cE_{\bar{\sigma}}(\pi)^*J\otimes 1)u_{\pi\otimes \pi^{\mathrm{op}}}(u_{\bar{\sigma}}\otimes 1)\cE_\sigma^* \\
&=  (J\cE_{\bar{\sigma}}(\pi)^*J\otimes 1)\kappa_\pi^*u_{\pi\otimes \pi^{\mathrm{op}}}   (R_\sigma^*u_{\bar{\sigma}}\otimes 1)\cE_\sigma^* \\
&= \gamma_R(\cE_{\bar{\sigma}}(\pi)\otimes 1)\gamma_R(\lambda_\pi^*)  (R_\sigma^*u_{\bar{\sigma}}\otimes 1)\cE_\sigma^*
= \gamma_R(\tilde{\bar{\sigma}}(\lambda_\pi^*))  (R_\sigma^*u_{\bar{\sigma}}\otimes 1)\cE_\sigma^*
\end{align*}
for $\pi\in \Irr(\cC_0)$.
Thus $(R_\sigma^*u_{\bar{\sigma}}\otimes 1)\cE_\sigma^*$ gives an equivalence of 
${}_{\tilde{\sigma}} L^2(\cB_0)$ and $L^2(\cB_0)_{\tilde{\bar{\sigma}}}$, and hence we have 
$\overline{\eta(\sigma,\cE)}=\eta(\bar{\sigma},\cE_{\bar{\sigma}})$.

We next show the equivalence
${}_{\tilde{\sigma}} L^2(\cB_0)$ and $L^2(\cB_0)_{\tilde{\bar{\sigma}}}$
for $\sigma$ with $d(\sigma)=\infty$.
If we examine the proof for $\sigma$ with $d(\sigma)<\infty$, we only have to show
\[
 \kappa_\pi^*(R_\sigma^*u_{\bar{\sigma}}\otimes 1)(\cE_\sigma(\pi)^*\otimes 1)u_{\pi\otimes \pi^{\mathrm{op}}}
=\gamma_R(\cE_{\bar{\sigma}}(\pi)\otimes 1)\gamma_R(\lambda_\pi^*)  (R_\sigma^*u_{\bar{\sigma}}\otimes 1).
\]
Note that we can see 
$R_\sigma^*u_{\bar{\sigma}}=\sum_k  J \bar{w}_kJ  R_{\xi_k}^*u_{\bar{\xi}_k}{w}_k^*$,
and the right hand side is a well-defined unitary even if $d(\sigma)=\infty$.
By substituting this expression of $R_\sigma^*u_{\bar{\sigma}}$, 
the above equality can be written as 
\begin{align*}
&\sum_{k,l}
(J \bar{w}_kJ\otimes 1)
\kappa_\pi^*( R_{\xi_k}^*u_{\bar{\xi}_k}\otimes 1)
(\cE_\sigma(\pi)^*_{k,l}\otimes 1)u_{\pi\otimes \pi^{\mathrm{op}}}({w}_l^*   \otimes 1) \\
&=
\sum_{k,l} 
( J \bar{w}_k\cE_{\bar{\sigma}}(\pi)_{k,l}^*J\otimes 1)
\kappa_\pi^*u_{\pi\otimes \pi^{\mathrm{op}}}
(  R_{\xi_l}^*u_{\bar{\xi}_l}w_l^*\otimes 1).
\end{align*}
Thus we only have to show
\[
\kappa_\pi^*( R_{\xi_k}^*u_{\bar{\xi}_k}\otimes 1)
(\cE_\sigma(\pi)^*_{k,l}\otimes 1)u_{\pi\otimes \pi^{\mathrm{op}}}
=
( J \cE_{\bar{\sigma}}(\pi)_{k,l}^*J\otimes 1)
\kappa_\pi^*u_{\pi\otimes \pi^{\mathrm{op}}}
(  R_{\xi_l}^*u_{\bar{\xi}_l}\otimes 1)
\]
for all $k,l$. This equation can be verified in a similar way as in $d(\sigma)<\infty$ case.
\hfill$\Box$

\medskip

\noindent
\textbf{Remark.} When $|\Irr(\cC)|<\infty$, $Z_{\cC_0}(\cC)$ is a rigid C$^*$-tensor category.
However in the case $|\Irr(\cC)|=\infty$, $Z_{\cC_0}(\cC)$ contains an element with infinite dimension.
(See Example \ref{exam:Zunitary}.)
So it is a tensor category, but it is 
not rigid. (Even in this case,  we can define a conjugate object as in  Proposition \ref{prop:eta-inter}.)

\medskip

Let us assume $|\Irr(\cC)|<\infty$. 
One of the main results of \cite{Kw-Dcenter} is
the existence of 1-to-1 correspondence
between $\mathrm{Irr} (Z_{\cC_0}(\cC))$ and the set of minimal central projections of the relative Tube algebra
for $\cC_0\subset \cC$.  
First, we will study this fact.

By definition of the relative tube algebra and $\cC_g\cap \cC_h=\emptyset$, $g\ne h$,
we have
\[
\mathrm{Tube}(\cC_0, \cC)=\bigoplus_{\substack{\xi\in \cC_0,  \pi, \rho\in \cC}}
(\pi\xi,\xi\rho)=
\bigoplus_{g\in G}\bigoplus_{\substack{\xi\in \cC_0,  \pi, \rho\in \cC_g}}
(\pi\xi,\xi\rho).
\]

Let
\[
\mathrm{Tube}(\cC_0, \cC)_g:= \bigoplus_{\substack{\xi\in \cC_0,  \pi, \rho\in \cC_g}}
(\pi\xi,\xi\rho).
\]
We can see that $\mathrm{Tube}(\cC_0,\cC)=\bigoplus_{g\in G}\mathrm{Tube}(\cC_0,\cC)_g$, and
$\mathrm{Tube}(\cC_0, \cC)_g$ is a (non-unital) $*$-subalgebra of $\mathrm{Tube}(\cC_0,\cC)$.

Let
\[
Z_{\cC_0}(\cC)_g:=
\left\{\theta\in Z_{\cC_0}(\cC)\mid \theta\sim \eta(\sigma,\mathcal{E}_\sigma) \mbox{ for some }\sigma\in
\cC_g \right\}.
\]
By the same argument in the proof of \cite[Lemma 4.5]{Iz-LR},
\cite[Theorem 3.2]{Kw-Dcenter}, we have
\begin{align*}
Z_{\cC_0}(\cC)_g& =
\left\{\theta\in Z_{\cC_0}(\cC)\mid \theta\iota \sim \iota(\sigma\otimes \id) \mbox{ for some }\sigma\in
\cC_g \right\} \\
&=\left\{\theta\in Z_{\cC_0}(\cC)\mid \theta\prec \iota(\sigma\otimes \id)\bar{\iota}\, \mbox{ for some }\sigma\in
\cC_g \right\} \,\,  (\mbox{when }d(\iota)<\infty ).
\end{align*}
It is clear that
\[
Z_{\cC_0}(\cC)=\bigoplus_{g\in G}Z_{\cC_0}(\cC)_g,\,\,
Z_{\cC_0}(\cC)_gZ_{\cC_0}(\cC)_h\subset Z_{\cC_0}(\cC)_{gh},\,\,
\overline{Z_{\cC_0}(\cC)_g}=Z_{\cC_0}(\cC)_{g^{-1}}.
\]

\begin{lem}\label{lem:diffgrade}
Suppose $|\Irr(\cC)|<\infty$. Then
there exists 1 to 1 correspondence between $\mathrm{Irr}(Z_{\cC_0}(\cC)_g)$ and the
set of minimal central projections of $\mathrm{Tube}(\cC_0,\cC)_g$.
\end{lem}
\textbf{Proof.}
This can be verified by examining the proof of \cite[Theorem 4.10]{Iz-LR}.
Here we present a proof based on the argument in \cite[Remark 5.1]{Mu-subII}.
(Also see Appendix \ref{sec:TubeReal}.)
By the characterization of $Z_{\cC_0}(\cC)_g$, all elements of $\mathrm{Irr}(Z_{\cC_0}(\cC)_g)$
are contained in $\bigoplus_{\pi\in \cC_g}\iota\left(\pi \otimes \id\right)\bar{\iota}$.
Hence the set of minimal central projections of
\[
\left(\bigoplus_{\pi\in \cC_g}\iota\left(\pi \otimes \id\right)\bar{\iota},
\bigoplus_{\rho\in \cC_g}\iota\left(\rho \otimes \id\right)\bar{\iota}
\right)
\]
corresponds to $\mathrm{Irr}(Z_{\cC_0}(\cC)_g)$.

By the Frobenius reciprocity, 
\begin{align*}
\left(\bigoplus_{\pi\in \cC_g}\iota(\pi\otimes \id)\bar{\iota}, \bigoplus_{\rho\in \cC_g}\iota(\rho\otimes \id)\bar{\iota}\right)&\cong
\bigoplus_{\pi,\rho\in \cC_g}\left((\pi\otimes \id)\bar{\iota}\iota,
\bar{\iota}\iota(\rho\otimes \id)\right) \\
&\cong \bigoplus_{\pi,\rho\in \cC_g, \xi,\eta\in \cC_0}
\left((\pi\otimes \id)(\xi\otimes \xi^{\mathrm{op}}),
(\eta\otimes \eta^{\mathrm{op}})(\rho\otimes \id)\right) \\
&\cong \bigoplus_{\pi,\rho\in \cC_g, \xi\in \cC_0}
\left(\pi\xi\otimes \xi^{\mathrm{op}},
\xi\rho\otimes \xi^{\mathrm{op}}\right) \\
&\cong\bigoplus_{\pi,\rho\in \cC_g, \xi\in \cC_0}(\pi\xi,\xi\rho)
\end{align*}
as Hilbert spaces.
If we carefully examine the definition of a tube algebra, we 
obtain the desired result.  (See \S \ref{sec:TubeReal} on the detail.)
\hfill$\Box$

\bigskip

Our next purpose is to show the existence of a $G$-braiding on $Z_{\cC_0}(\cC)$.
Recall that 
there exists an action $\beta$ of $G$ on $\cB_0$
arising from $\cC$ as explained in \S \ref{sec:LoRe}. 
We will prepare several lemmas to construct a $G$-braiding. 

\begin{lem}\label{lem:relationgenerator2}
For 
$\pi\in \cC_g$, define
$\hlambda_\pi\in \cB_0$ by $\gamma_L(\hlambda_\pi)=\gamma_L(\lambda_\pi)U_g^*$.
(Note $\gamma_L(\lambda_\pi)U_g^*\in B(\fH_0)\cap\gamma_R(\cB_0)'=\gamma_L(\cB_0)$.)
Then
$\hlambda_\pi\in (\beta_g\iota,\iota (\pi\otimes \pi^\mathrm{op}))$.
For $\pi,\pi'\in \cC_g$, $\rho\in \cC_h$, we have 
\[
\hlambda_\pi\beta_g(\hlambda_\rho)
=\hlambda_{\pi\rho}, \,\,
\hlambda_\pi^*=\beta_g(R_{\pi\otimes \pi^\mathrm{op}}^*\hlambda_{\bar{\pi}}),\,\,\, 
(T\otimes 1)\hlambda_\pi=(1\otimes j(T^*))\hlambda_{\pi'},\,\, T\in (\pi,\pi').
\]
In particular, we have 
$\hlambda_\pi\beta_g(\hlambda_\rho)=
\sum_{\xi\in \cC_{gh}}\tilde{T}_{\pi,\rho}^\xi \hlambda_\xi$, $\pi\in \Irr(\cC_g)$, $\rho\in \Irr(\cC_h)$, 
and $\hlambda_\pi^*\hlambda_\pi=d(\pi)$.
\end{lem}
\textbf{Proof.}
The facts $\hlambda_\pi\in (\beta_g\iota,\iota (\pi\otimes \pi^\mathrm{op}))$ 
and $(T\otimes 1)\hlambda_\pi=(1\otimes j(T^*))\hlambda_{\pi'}$ can be easily verified.
The other properties can be shown as follows;
\[
 \gamma_L(\hlambda_\pi)\gamma_L(\beta_g(\hlambda_\rho))=
 \gamma_L(\lambda_\pi)\gamma_L(\hlambda_\rho)U_g^*=\gamma_L(\lambda_{\pi\rho})U_{gh}^*
=\gamma_L(\hlambda_{\pi\rho}),
\]
\[
 \gamma_L(\hlambda_\pi^*)=U_g\gamma_L(\lambda^*_\pi)=
U_g\gamma_L((R_\pi^*\otimes j(R_\pi^*))\lambda_{\bar{\pi}})U_{g^{-1}}^*U_{g}^*=
\gamma_L(\beta_g((R_\pi^*\otimes j(R_\pi^*))\hlambda_{\bar{\pi}})).
\]

\hfill$\Box$

\begin{lem}\label{lem:conjugate}
For $\pi\in \mathrm{Irr}(\cC_g)$, $\iota(\pi\otimes \id)$ is irreducible, and
$\beta_g\iota(\id\otimes \bar{\pi}^{\mathrm{op}})\sim 
\iota(\pi\otimes \id)$.
The latter equivalence is given by $(1\otimes j(\bar{R}_\pi^*))\hlambda_\pi $.
\end{lem}
\textbf{Proof.}
Let $a=\sum_{\rho\in \cC_0}a_\rho \lambda_\rho\in (\iota(\pi\otimes \id),
\iota(\pi\otimes \id))$.
Then we have 
\[
 \sum_{\rho\in \cC_0}(\pi\otimes \id)(x)a_\rho \lambda_\rho=
 \sum_{\rho\in \cC_0}a_\rho \lambda_\rho(\pi\otimes \id)(x)=
 \sum_{\rho\in \cC_0}a_\rho(\rho\pi\otimes \rho^{\mathrm{op}})(x)\lambda_\rho.
\]
Comparing coefficients, we get 
 $a_\rho\in (\rho\pi\otimes \rho^{\mathrm{op}}, \pi\otimes \id)$ for all $\rho\in \Irr(\cC_0)$.
Hence $a_\rho=0$ for $\rho\ne \id$, and $a_{\id}\in (\pi\otimes \id,\pi\otimes \id)=\mathbb{C}$.
We can see that $\left((1\otimes j(\bar{R}_\pi^*))\hlambda_\pi\right)^*= 
\beta_g((R_\pi^*\otimes 1)\hlambda_{\bar{\pi}})$, and  
$(1\otimes j(\bar{R}_\pi^*))\hlambda_\pi $ is a unitary in 
$(\beta_g\iota(\id\otimes \bar{\pi}^{\mathrm{op}}),\iota(\pi\otimes \id) )$ 
by Lemma \ref{lem:relationgenerator2}. 
\hfill$\Box$

\medskip

\begin{prop}\label{prop:Vunitary}
Let $V_\sigma:=\left(1\otimes j(\bar{R}_\sigma^*)\right)\hlambda_\sigma$, $\sigma\in \cC$.
Then $V_\sigma$ is a unitary
in $(\beta_g(\id\otimes \bar{\sigma}^{\mathrm{op}}),\sigma\otimes \id)$ for $\sigma\in \cC_g$. 
Moreover we have \\
$(\mathrm{i})$ For $\sigma\in \cC_g$ and $\rho\in \cC$,
 $V_\sigma\beta_g(V_\rho)=V_{\sigma\rho}$, where we choose $\bar{R}_{\sigma\rho}$
as $\sigma(\bar{R}_\rho)\bar{R}_\sigma$. \\
$(\mathrm{ii})$ 
For $T\in (\sigma,\sigma')$, define $\hat{T}\in (\bar{\sigma},\bar{\sigma'})$ by
$\hat{T}:=
R_\sigma^*\bar{\sigma}(T^*)\bar{\sigma}(\bar{R}_{\sigma'})=
\bar{\sigma'}(\bar{R}_\sigma)^*\bar{\sigma'}(T^*)R_{\sigma'}$. \\
Then we have $(T\otimes 1)V_\sigma=V_{\sigma'}\beta_g\left(1\otimes j(\hat{T})\right)$ and
$(T^*\otimes 1)V_{\sigma'}=V_\sigma\beta_g\left(1\otimes j(\hat{T}^*)\right)$. \\
$(\mathrm{iii})$ $\tilde{\sigma}=\Ad V_\sigma \circ \beta_g \tilde{\bar{\sigma}}^{\mathrm{op}}$ for 
$\tilde{\sigma}\in Z_{\cC_0}(\cC)_g$. 
\end{prop}
\textbf{Proof.}
At first we remark that
\[
V_\sigma=\sum_{k}(w_k\otimes 1 )V_{\xi_k}\beta_g(1\otimes j(\bar{w}_k^*))
\]
for irreducible decomposition $\sigma(x)=\sum_kw_k\xi_k(x)w_k^*$ and 
$\bar{\sigma}(x)=\sum_k\bar{w}_k\bar{\xi}_k(x)\bar{w}_k^*$ 
with $d(\sigma)<\infty$.
So for $d(\sigma)=\infty$, $V_\sigma$ should be understood by the above formula.
The proof presented below can be easily modified for $\sigma$ with 
$d(\sigma)=\infty$ as in Proposition \ref{prop:eta-inter}.

By Proposition \ref{lem:relationgenerator2} (or Lemma \ref{lem:conjugate}),  
we can see that $V_\sigma$ is a unitary in $(\beta_g(\id\otimes \bar{\sigma}^{\mathrm{op}}),\sigma\otimes \id)$. 
Indeed, $V_\sigma^*$ is given by
$V_\sigma^*=\beta_g((R_\sigma^*\otimes 1)\hlambda_{\bar{\sigma}})$.

(i) By Proposition \ref{lem:relationgenerator2}, 
\begin{align*}
V_\sigma\beta_g(V_\rho)&=
(1\otimes j(\bar{R}_\sigma^*))\hlambda_\sigma
\beta_g\left((1\otimes j(\bar{R}_\rho^*))\hlambda_\rho\right)
=
(1\otimes j(\bar{R}_\sigma^*\sigma(\bar{R}_\rho)^*)\hlambda_\sigma
\beta_g\left(\hlambda_\rho\right) \\
&=
(1\otimes j(\bar{R}_\sigma^*\sigma(\bar{R}_\rho)^*)\hlambda_{\sigma\rho}\\
&=V_{\sigma\rho}.
\end{align*}

(ii)  Note that 
$\sigma'(\hat{T}^*)\bar{R}_{\sigma'}=T\bar{R}_{\sigma}$ 
holds. 
Combining with Proposition \ref{lem:relationgenerator2}, we get 
\begin{align*}
(T\otimes 1)V_\sigma&=\left(T\otimes j(\bar{R}_\sigma^*) \right)\hlambda_\sigma=
\left(1\otimes j(\bar{R}_\sigma^*T^*) \right)\hlambda_{\sigma'}\\
&=\left(1\otimes j(\bar{R}_{\sigma'}^*\sigma'(\hat{T})) \right)\hlambda_{\sigma'}=
\left(1\otimes j(\bar{R}_{\sigma'}^*) \right)\hlambda_{\sigma'}\beta_g\left(1\otimes j(\hat{T})\right)\\
&=V_{\sigma'}\beta_g\left(1\otimes j(\hat{T})\right).
\end{align*}
In a similar way, we have $(T^*\otimes 1)V_{\sigma'}=V_\sigma\beta_g\left(1\otimes j(\hat{T}^*)\right)$
by using $\sigma(\hat{T})\bar{R}_\sigma=T^*\bar{R}_{\sigma'}$.

(iii) It suffices to show
$  V_\sigma\beta_g\tilde{\bar{\sigma}}^{\mathrm{op}}(\lambda_\rho)
= \tilde{\sigma}(\lambda_\rho)V_\sigma
$ for $\rho\in \mathrm{Irr}(\cC_0)$.
Since $\hlambda_\rho=\lambda_\rho$, $\rho\in \cC_0$, we have
\begin{align*}
 V_\sigma\beta_g\tilde{\bar{\sigma}}^{\mathrm{op}}(\lambda_\rho)
&=
(1\otimes j(\bar{R}_\sigma^*))\hlambda_\sigma
\beta_g ((1\otimes j (\cE_{\bar{\sigma}}(\rho)^*))\lambda_\rho) \\
&=
(1\otimes j(\bar{R}_\sigma^* \sigma(\cE_{\bar{\sigma}}(\rho)^*)))
\hlambda_\sigma\beta_g (\lambda_\rho) 
=(1\otimes j(\rho(\bar{R}_\sigma^*) \cE_{\sigma}(\rho)))
\hlambda_{\sigma\rho} \\
&=
(\cE_{\sigma}(\rho)^* \otimes j(\rho(\bar{R}_\sigma^*) ))
\hlambda_{\rho\sigma} \,\,\,  
=
(\cE_{\sigma}(\rho)^* \otimes 1)\lambda_\rho
(1 \otimes j(\bar{R}_\sigma^*) )
\hlambda_{\sigma} \\
&=\tilde{\sigma}(\lambda_\rho)V_\sigma
\end{align*}
 by Proposition \ref{lem:relationgenerator2}.  
We comment that  (iii)  is the generalization of \cite[Theorem 6.4]{Iz-LR}.
(Also see the comment after \cite[Proposition 2.6]{BEK-LR}.)
\hfill $\Box$

\medskip

\noindent
\textbf{Remark.} The definition of $V_\sigma$ depends on the choice of $(R_\sigma,\bar{R}_\sigma)$.

\medskip

Proposition  \ref{prop:Vunitary}
suggests how to construct a $G$-braiding operator, since we have
\[
\beta_g\eta(\rho,\mathcal{E}_\rho)\beta^{-1}_g\eta(\sigma,\mathcal{E}_\sigma)\sim
\beta_g\eta(\rho,\mathcal{E}_\rho){\eta^{\mathrm{op}}(\bar{\sigma}, \mathcal{E}_{\bar{\sigma}})}
=\beta_g{\eta^{\mathrm{op}}(\bar{\sigma}, \mathcal{E}_{\bar{\sigma}})}
\eta(\rho,\mathcal{E}_\rho)\sim
\eta(\sigma,\mathcal{E}_\sigma)
\eta(\rho,\mathcal{E}_\rho).
\]
We observe that our construction of a $G$-braiding is similar to the construction of braiding in AQFT.

\begin{thm}\label{thm:Gbraid0}
Let us denote $\beta_g\tilde{\sigma}\beta_g^{-1}$ by $g[\tilde{\sigma}]$.
Define
\[
\mathcal{E}(\tilde{\sigma},\tilde{\rho}):=g[\tilde{\rho}](V_\sigma)V_\sigma^*, \,\,
\tilde{\sigma}\in Z_{\cC_0}(\cC)_g, \,\,
\tilde{\rho}\in Z_{\cC_0}(\cC).
\]
and extend it to $\mathcal{E}(\theta_1,\theta_2)$ canonically for all $\theta_1,\theta_2\in Z_{\cC_0}(\cC)$.
Then \\
$(1)$ $\cE(\tilde{\sigma},\tilde{\rho})$ does not depend on the choice of $(R_\sigma,\bar{R}_\sigma)$. \\
$(2)$ $\cE(\theta_1,\theta_2)\in (\theta_1\theta_2,g[\theta_2]\theta_1)$ for 
$\theta_1\in Z_{\cC_0}(\cC)_g$, $\theta_2\in Z_{\cC_0}(\cC)$. \\
$(3)$ $g[\theta]\in Z_{\cC_0}(\cC)_{ghg^{-1}}$ for $\theta\in Z_{\cC_0}(\cC)_h$.
\end{thm}
\textbf{Proof.}
(1)  Let $(R'_\sigma, \bar{R}'_\sigma)$ be another choice of a solution of conjugate equation for $\sigma,\bar{\sigma}$.
Since there exists a unitary $u\in (\bar{\sigma},\bar{\sigma})$ such that $R_\sigma'=uR_\sigma$, $\bar{R}_\sigma'=\sigma(u)\bar{R}_\sigma$,
we can see that $\mathcal{E}(\tilde{\sigma},\tilde{\rho})$ does not depend on the choice of $\bar{R}_\sigma$.
(Note that the definition of $\cE(\tilde{\sigma},\tilde{\rho})$ 
depends only on $\sigma$.)

(2) Since $\tilde{\sigma}^{\mathrm{op}}$ and $\tilde{\rho}$ commute, we have
\begin{align*}
g[\tilde{\rho}]\tilde{\sigma}(a)\mathcal{E}(\tilde{\sigma},\tilde{\rho})&=
\beta_g\tilde{\rho}\beta_g^{-1}\tilde{\sigma}(a)\beta_g\tilde{\rho}\beta_g^{-1}(V_\sigma)V_\sigma^*
=\beta_g\tilde{\rho}\beta_g^{-1}(V_\sigma)\beta_g\tilde{\rho}\tilde{\bar{\sigma}}^{\mathrm{op}}(a)
V_\sigma^* \\
&=\beta_g\tilde{\rho}\beta_g^{-1}(V_\sigma)\beta_g\tilde{\bar{\sigma}}^{\mathrm{op}}\tilde{\rho}(a)V_\sigma^*
=\beta_g\tilde{\rho}\beta_g^{-1}(V_\sigma)V_\sigma^* \tilde{\sigma}\tilde{\rho}(a) \\
&= \mathcal{E}(\tilde{\sigma},\tilde{\rho})\tilde{\sigma}\tilde{\rho}(a)
\end{align*}
for all $a\in \cB_0$ by Proposition \ref{prop:Vunitary}.
Hence $\mathcal{E}(\tilde{\sigma},\tilde{\rho}) \in (\tilde{\sigma}\tilde{\rho}, g[\tilde{\rho}]\tilde{\sigma})$.

(3) Let  $\tilde{\sigma}\in Z_{\cC_0}(\cC)_g$,  $\tilde{\rho}\in Z_{\cC_0}(\cC)_h$,
and take $\pi\in \cC_g$ with 
$\sigma \succ \pi$, $d(\pi)<\infty$.
Then $\tilde{\sigma}\tilde{\rho}\iota \sim g[\tilde{\rho}]\tilde{\sigma}\iota$ implies
\[
\iota(\sigma\rho\otimes \id)\sim g[\tilde{\rho}]\iota ({\sigma}\otimes \id)\succ
g[\tilde{\rho}]\iota ({\pi}\otimes \id).
\]
Thus we get
\[
\iota(\sigma\rho\bar{\pi}\otimes \id)\succ 
g[\tilde{\rho}]\iota ({\pi}\bar{\pi}\otimes \id)\succ g[\tilde{\rho}]\iota .
\]
This implies $g[\tilde{\rho}]\iota \sim \iota(\xi\otimes \id)$ for some $\xi\in \cC_{ghg^{-1}}$,
and $g[\tilde{\rho}]\in  Z_{\cC_0}(\cC)_{ghg^{-1}}$.
\hfill$\Box$

\medskip

\noindent
\textbf{Remark.} When $[\cB_0:\cA]<\infty$, we can use the Frobenius reciprocity for $\iota$, $\bar{\iota}$.
Then Theorem \ref{thm:Gbraid0} (3) can be verified in a similar way as in 
\cite[Theorem 3.2]{Kw-Dcenter}.

\medskip

Any element of $Z_{\cC_0}(\cC)$ is of the form $\eta(\sigma, \mathcal{E}_\sigma)$ up to unitary equivalence.
Hence Theorem \ref{thm:Gbraid0} (3)
implies
$\beta_k\eta(\sigma,\mathcal{E})\beta^{-1}_k$ is equivalent to
$\eta(\sigma^{(k)},\mathcal{E}_{\sigma^{(k)}})$ for some $\sigma^{(k)}$ and its half braiding $\mathcal{E}_{\sigma^{(k)}}$. However, we do not know the explicit expression of $\sigma^{(k)}$ and $\mathcal{E}_{\sigma^{(k)}}$
in terms of $\sigma$, $\mathcal{E}_\sigma$, $k$. 
Any way, we know there exists a unitary in $(\beta_k\tilde{\sigma}\beta_k^{-1},\tilde{\sigma}^{(k)})$. 
To describe the intertwiner space $(\beta_k\tilde{\sigma}\beta_k^{-1},\tilde{\sigma}^{(k)})$, 
we first prepare the following Lemma, whose proof is similar to that of Lemma \ref{lem:halfglobal}.

\begin{lem}\label{lem:actinInter}
Let $\{Z_\sigma(\pi)\}_{\pi\in \mathrm{Irr}(\cC_k)}\subset \cM$ be a family of unitaries satisfying
the following conditions;
\begin{align*}
\mathrm{(A1)} &\,\,  Z_\sigma(\pi)\in (\sigma^{(k)}\pi,\pi \sigma) ,\,\, \pi \in \mathrm{Irr}(\cC_k),  \\
\mathrm{(A2)} &\,\,  TZ_\sigma(\rho)=\pi(\mathcal{E}_{\sigma}(\xi))Z_\sigma(\pi)\sigma^{(k)}(T), \,\, 
\xi\in \cC_0, \,\, \pi,\rho\in \cC_k, \,\,
T\in (\rho,\pi\xi), \\
\mathrm{(A3)} & \,\, Z_\sigma(\rho)=\xi(Z_\sigma(\pi))\mathcal{E}_{\sigma^{(k)}}(\xi)\sigma^{(k)}(S), \,\,
\xi\in \cC_0, \,\, \pi,\rho\in \cC_k, \,\,
S\in (\rho,\xi\pi). 
\end{align*}
Define a unitary $Z_\sigma\in \End(\fH_{k,\cA})$  by $(Z_\sigma v)(\xi)=Z_\sigma(\xi)v(\zeta)$.
Then $Z_\sigma$ satisfies the following;
\begin{align*}
\mathrm{(B1)}& \,\, Z_\sigma \sigma^{(k)}\gamma_L(a)=\gamma_L(\sigma\otimes \id(a))Z_\sigma,\,\,\, a\in \cA \\
\mathrm{ (B2)}&\,\, \gamma_L(\lambda_\pi) Z_\sigma=\gamma_L(\cE_\sigma(\pi)\otimes 1)Z_\sigma \sigma^{(k)}\gamma_L(\lambda_\pi),\,\, \,\pi\in \cC_0, \\
\mathrm{ (B3)}&\,\, \kappa_\pi Z_\sigma=\pi(Z_\sigma)(\cE_{\sigma^{(k)}}(\pi)\otimes 1)\sigma^{(k)}(\kappa_\pi),\,\, \,\pi\in \cC_0.
\end{align*}
Conversely, if 
$Z_\sigma\in \End(\fH_{k,\cA})$  
satisfies (B), then $Z_\sigma $ comes from 
$\{Z_\sigma(\pi)\}$ satisfying (A).
\end{lem}

We can now describe an intertwiner space
$\left(\eta(\sigma^{(k)},\mathcal{E}_{\sigma^{(k)}}),
\beta_k\eta(\sigma,\mathcal{E}_\sigma)\beta_k^{-1}
\right)$.

\begin{thm}\label{thm:actionInter}
Let $Z_\sigma$ be as in Lemma \ref{lem:actinInter}. Then we have the following \\
$(1)$ $U_kZ_\sigma \sigma^{(k)} (U_k^*)\mathcal{E}_{\sigma^{(k)}}^*\in \gamma_R(\cB_0)'=\gamma_L(\cB_0)$. \\
$(2)$ Let $u\in \cB_0$ be a unitary with $\gamma_L(u)=
U_kZ_\sigma \sigma^{(k)}(U_k^*)\mathcal{E}_{\sigma^{(k)}}^*$.
Then
$\beta_k\tilde{\sigma}\beta_k^{-1}=\Ad u\circ \tilde{\sigma}^{(k)}$. \\
$(3)$ If a unitary $u\in \cB_0$ satisfies
$\beta_k\tilde{\sigma}\beta_k^{-1}=\Ad u\circ \tilde{\sigma}^{(k)}$,
then
$Z:=U_k^*\gamma_L(u)\mathcal{E}_{\sigma^{(k)}}\sigma^{(k)}(U_k)$ satisfies the condition (2)
in Lemma \ref{lem:actinInter}.
\end{thm}
\textbf{Proof.}
(1) It is clear that $U_kZ_\sigma\sigma^{(k)}(U_k^*)\cE_{\sigma^{(k)}}^*\in \gamma_R(\cA)'$.
By Lemma \ref{lem:halfglobal} and Lemma \ref{lem:actinInter} (B3), 
\begin{align*}
&\gamma_R(\lambda_\pi)^*U_kZ_\sigma\sigma^{(k)}(U_k^*)\cE_{\sigma^{(k)}}^*
=
U_k\gamma_R(\lambda_\pi)^*Z_\sigma\sigma^{(k)}(U_k^*)\cE_{\sigma^{(k)}}^* \\
&=U_k\kappa^*_\pi u_{\pi\otimes \pi^{\mathrm{op}}} Z_\sigma\sigma^{(k)}(U_k^*)\cE_{\sigma^{(k)}}^* 
=U_k\kappa^*_\pi \pi(Z_\sigma)u_{\pi\otimes \pi^{\mathrm{op}}}\sigma^{(k)}(U_k^*)\cE_{\sigma^{(k)}}^* \\
&=
U_kZ_\sigma\sigma^{(k)}(\kappa_\pi^*)(\cE_{\sigma^{(k)}}(\pi)^*\otimes 1)u_{\pi\otimes \pi^{\mathrm{op}}}\sigma^{(k)}(U_k)^*
\cE_{\sigma^{(k)}}^*   \\
&=
U_kZ_\sigma\sigma^{(k)}(\kappa_\pi^*)(\cE_{\sigma^{(k)}}(\pi)^*\otimes  1)\pi\sigma^{(k)}(U_k)^*u_{\pi\otimes \pi^{\mathrm{op}}}
\cE_{\sigma^{(k)}}^* \\ 
&=
U_kZ_\sigma\sigma^{(k)}(\kappa_\pi^*)\sigma^{(k)}\pi(U_k)^*(\cE_{\sigma^{(k)}}(\pi)^*\otimes 1)
u_{\pi\otimes \pi^{\mathrm{op}}} \cE_{\sigma^{(k)}}^*  \\
&=
U_kZ_\sigma\sigma^{(k)}(\kappa_\pi^*\pi(U_k)^*)(\cE_{\sigma^{(k)}}(\pi)^*\otimes 1)
u_{\pi\otimes \pi^{\mathrm{op}}} \cE_{\sigma^{(k)}}^* \\ 
&=
U_kZ_\sigma\sigma^{(k)}(U_k^*\kappa_\pi^*)(\cE_{\sigma^{(k)}}(\pi)^*\otimes 1)
\pi( \cE_{\sigma^{(k)}})^*u_{\pi\otimes \pi^{\mathrm{op}}}  
=
U_kZ_\sigma\sigma^{(k)}(U_k^*)
 \cE_{\sigma^{(k)}}^*\gamma_R(\lambda_\pi^*).
\end{align*}
Hence we obtain the conclusion. \\
\noindent
(2)  Lemma \ref{lem:actinInter} (B1) and (B2) yield
$Z_\sigma\sigma^{(k)}\gamma_L(a)Z_\sigma^*=\gamma_L(\tilde{\sigma}(a))$
for $a\in \cB_0$.  Then 
\begin{align*}
& U_kZ_\sigma\sigma^{(k)}(U_k^*)\cE_{\sigma^{(k)}}^*
\gamma_L(\tilde{\sigma}^{(k)}\beta_k(a))
=
 U_kZ_\sigma\sigma^{(k)}(U_k^*)\sigma^{(k)}\left(\gamma_L(\beta_k(a))\right)\cE_{\sigma^{(k)}}^* \\
&=
 U_kZ_\sigma\sigma^{(k)}(\gamma_L(a))\sigma^{(k)}(U_k^*)
\cE_{\sigma^{(k)}}^*
=
 U_k\gamma_L(\tilde{\sigma}(a))Z_\sigma\sigma^{(k)}(U_k^*)
\cE_{\sigma^{(k)}}^*\\
&=\gamma_L(\beta_k\tilde{\sigma}(a))
U_kZ_\sigma\sigma^{(k)}(U_k^*)\cE_{\sigma^{(k)}}^*.
\end{align*}
Thus we get $\beta_k\tilde{\sigma}\beta_k^{-1}=\Ad u\circ \tilde{\sigma}^{(k)}$. 

If we examine the above proof, we can easily see that the statement (3) holds.
\hfill$\Box$

\medskip

\noindent
\textbf{Remark} When $k=0$, we can assume $u=1$ in Theorem \ref{thm:actionInter}.
In this case, we have $Z_\sigma(\pi)=\cE_\sigma(\pi)$.

\begin{exam}\label{exam:Zunitary}
\upshape 
We present an example of 
Theorem \ref{thm:actionInter}.
Let $\{s_\xi^k\}_{\xi\in \cC_k}\subset \cM$ be an set of isometries with 
$\sum_{\xi\in \cC_k}s_\xi^ks_\xi^{k*}=1$, $k\in G$. Fix $\mu\in \cC_g$, and 
define $\theta_\mu^{(k)}\in \End(\cM)$, and a unitary $Z_{\theta_\mu^{(k)}}^h(\pi)\in \cM$, 
$\pi\in \cC_h$, by 
\begin{align*}
 \theta_\mu^{(k)}(x):=\sum_{\xi\in \cC_k}s_\xi^k\xi \mu\bar{\xi}(x)s_\xi^{k*}, \,\,\,
Z_{\theta_\mu^{(k)}}^h(\pi):=
\sum_{\xi\in \cC_k,\zeta\in \cC_{hk}}\pi(s_\xi^k)T_{\pi,\xi}^{\zeta,e}
\zeta\mu(T_{\bar{\zeta},\pi}^{\bar{\xi},e})^*s_\zeta^{hk*},
\end{align*}
where $\{T_{\pi,\xi}^{\zeta,e}\}_e\subset (\zeta,\pi\xi)$  and $\{T_{\bar{\zeta},\pi}^{\bar{\xi},e}\}_e
\subset (\bar{\xi},\bar{\zeta}\pi)$ 
are both CONS, which are related by the Frobenius reciprocity, i.e.,
\[
 T_{\bar{\zeta},\pi}^{\bar{\xi},e}=\sqrt{\frac{d(\xi)}{d(\zeta)}}\bar{\zeta}\pi(\bar{R}_\xi^*)T_{\pi,\xi}^{\zeta,e}R_\zeta.
\]

We can see $\theta_\mu^{(k)}\in \cC_{kgk^{-1}}$, 
$Z_{\theta_\mu^{(k)}}^{h}(\pi)\in (\theta_\mu^{(hk)}\pi,\pi\theta_\mu^{(k)})$, and 
\[
\pi(Z^{h_2}_{\theta_{\mu}^{(k)}}(\rho))Z_{\theta_{\mu}^{(h_2k)}}^{h_1}(\pi)
=Z_{\theta_{\mu}^{(k)}}^{h_1h_2}(\pi\rho),\,\,\,
 TZ_{\theta_{\mu}^{(k)}}^{h}(\sigma)=Z_{\theta_{\mu}^{(k)}}^{h}(\sigma')\theta_\mu^{(hk)}(T)
\]
for 
$\pi\in \cC_{h_1}$, $\rho\in \cC_{h_2}$, $\sigma,\sigma'\in \cC_h$,
$T\in (\sigma,\sigma')$, e.g., \cite[\S 3.2]{Nesh-Yama-DrinRep}.
Note $\theta_\mu^{(k)}$ is of infinite index if $|\Irr(\cC)|=\infty$.

In particular, $\{Z_{\theta_{\mu}^{(k)}}^0(\pi)\}_{\pi\in \cC_0}$ is a half braiding for $\theta_\mu^{(k)}$, and 
$\{Z_{\theta_{\mu}^{(0)}}^h(\pi)\}_{\pi\in \cC_h}$ satisfies the condition (A) of  Lemma \ref{lem:actinInter}.
By Theorem \ref{thm:actionInter}, $\beta_k\tilde{\theta}_\mu^{(0)}\beta_k^{-1}\sim \tilde{\theta}_\mu^{(k)}$.
We remark that $\tilde{\theta}_\mu^{(0)}\sim \iota(\mu\otimes \id) \bar{\iota}$, cf,  \cite[Proposition 6.2]{Iz-LR}.
(See Appendix.) 
\end{exam}

\medskip

Let $\{Z_\sigma(\pi)\}_{\pi\in \cC_k}$ be as in Lemma \ref{lem:actinInter} (1), and define 
a unitary $Z_\sigma^{\mathrm{op}}\in B(\fH_k)$ by
$\left( Z_\sigma^{\mathrm{op}} v\right)(\pi)=\big(1\otimes j(Z_\sigma(\pi))\big)v(\pi)$.
In a similar way, we can construct a unitary $w\in \cB_0$ such that
$\beta_k\tilde{\sigma}^{\mathrm{op}}\beta_k^{-1}
=\Ad(w)\circ \tilde{\sigma}^{(k),\mathrm{op}}$
by $\gamma_L(w)=U_kZ_\sigma^{\mathrm{op}}\sigma^{(k),\mathrm{op}}(U_k^*)
\mathcal{E}_{\sigma^{(k)}}^{\mathrm{op}*}$.

The following Lemma will be used in the proof of Theorem \ref{thm:Gbraid1}.
\begin{lem}\label{lem:Intercommute}
Let
$v\in (\tilde{\rho}^{(k)},k[\tilde{\rho}])$ and
$w\in (\tilde{\sigma}^{(k),\mathrm{op}}, k[\tilde{\sigma}^{\mathrm{op}}])$
be unitaries constructed in Theorem \ref{thm:actionInter}.
Then
we have $v\tilde{\rho}^{(k)}(w)=
w\tilde{\sigma}^{(k),\mathrm{op}}(v)$.
\end{lem}
\textbf{Proof.} Let $\{Z_\rho(\pi)\}_{\pi\in \cC_k}$ and $\{Z_\sigma(\pi)\}_{\pi\in \cC_k}$ be
as in Lemma \ref{thm:actionInter} (1) for $\rho$ and $\sigma$ respectively.
By the construction in Theorem \ref{thm:actionInter},
\[
\gamma_L(v)=U_kZ_\rho(\rho^{(k)}\otimes \id)(U_k^*)\mathcal{E}_{\rho^{(k)}}^*, \,\,\,
\gamma_L(w)=U_kZ_\sigma^{\mathrm{op}}(\id\otimes \sigma^{(k),\mathrm{op}})(U_k^*)
\mathcal{E}_{\sigma^{(k)}}^{\mathrm{op}*}.
\]
By simple computation, we obtain
\begin{align*}
& \gamma_L(v)\gamma_L(\tilde{\rho}^{(k)}(w))=
\gamma_L(w)\gamma_L(\tilde{\sigma}^{(k),\mathrm{op}}(v))
=U_kZ_\rho Z_\sigma^{\mathrm{op}}(\rho^{(k)}\otimes \sigma^{(k),\mathrm{op}})(U_k^*)\cE_{\rho^{(k)}}^*
\mathcal{E}_{\sigma^{(k)}}^{\mathrm{op}*}.
\end{align*}

\hfill$\Box$

\medskip 

\begin{thm}\label{thm:Gbraid1}
A family of unitaries 
$\left\{\mathcal{E}(\tilde{\sigma},\tilde{\rho})\right\}_{\tilde{\sigma},\tilde{\rho}\in Z_{\cC_0}(\cC)}\subset \cB_0$ 
constructed in Theorem \ref{thm:Gbraid0} gives a 
$G$-braiding on $Z_{\cC_0}(\cC)$. Namely, we have 
\begin{align*}
\mathrm{(BF0)}\,\,\, &  \mathcal{E}(\tilde{\sigma},\tilde{\rho})\in (\tilde{\sigma}\tilde{\rho},g[\tilde{\rho}]\tilde{\sigma}),
\,\, \tilde{\sigma}\in Z_{\cC_0}(\cC)_g, \,\,\tilde{\rho}\in Z_{\cC_0}(\cC) ,\\ 
 \mathrm{(BF1)}\,\,\,  & 
g[\tilde{\xi}]\left(\mathcal{E}(\tilde{\sigma},\tilde{\eta})\right)\mathcal{E}(\tilde{\sigma},\tilde{\xi})
=\mathcal{E}(\tilde{\sigma},\tilde{\xi}\tilde{\eta}), \,\, \tilde{\sigma}\in Z_{\cC_0}(\cC)_g,  \\
&\beta_g(T)\mathcal{E}(\tilde{\sigma}, \tilde{\rho})=\mathcal{E}(\tilde{\sigma}, \tilde{\rho}')\tilde{\sigma}(T),\,\,
\tilde{\sigma}\in Z_{\cC_0}(\cC)_g,\,\, T\in (\tilde{\rho}, \tilde{\rho}'), \\
 \mathrm{(BF2)}\,\,\,  &
\mathcal{E}(\tilde{\sigma},h[\tilde{\xi}])\tilde{\sigma}\left(\mathcal{E}(\tilde{\rho},\tilde{\xi})\right)
=\mathcal{E}(\tilde{\sigma}\tilde{\rho},\tilde{\xi}), \,\,
\tilde{\sigma}\in Z_{\cC_0}(\cC)_g,\,\,
\tilde{\rho}\in Z_{\cC_0}(\cC)_h,\\ 
&g[\tilde{\xi}](T)\mathcal{E}(\tilde{\sigma}, \tilde{\xi})=
\mathcal{E}(\tilde{\sigma}', \tilde{\xi}) T, \,\, 
\tilde{\sigma},\tilde{\sigma}'\in Z_{\cC_0}(\cC)_g,\,\, T\in (\tilde{\sigma},\tilde{\sigma}'),\\
 \mathrm{(BF3)}\,\,\,  & \mathcal{E}(k[\tilde{\sigma}],k[\tilde{\rho}])=
\beta_k(\mathcal{E}(\tilde{\sigma},\tilde{\rho})).
\end{align*}

\end{thm}
\textbf{Proof.} The property (BF0) is already shown in Theorem \ref{thm:Gbraid0}.

(BF1).
The first relation of (BF1) can be verified as follows;
\begin{align*}
g[\tilde{\xi}]\left(\mathcal{E}(\tilde{\sigma},\tilde{\eta})\right)\mathcal{E}(\tilde{\sigma},\tilde{\xi})&=
\beta_g\tilde{\xi}\beta_g^{-1}\left(\beta_g\tilde{\eta}\beta_g^{-1}(V_\sigma)V_\sigma^*\right)
\beta_g\tilde{\xi}\beta_g^{-1}(V_\sigma)V_\sigma^*\\
&=\beta_g\tilde{\xi}\tilde{\eta}\beta_g^{-1}(V_\sigma)\beta_g\tilde{\xi}\beta_g^{-1}(V_\sigma^*)
\beta_g\tilde{\xi}\beta_g^{-1}(V_\sigma)V_\sigma^*\\
&=\mathcal{E}(\tilde{\sigma},\widetilde{\xi}\tilde{\eta}).
\end{align*}

Let $T\in (\tilde{\rho}, \tilde{\rho}')$.
Since $T\in \cM\otimes \mathbb{C}$, we have $\tilde{\bar{\sigma}}^{\mathrm{op}}(T)=T$. Thus
\begin{align*}
\beta_g(T)\mathcal{E}(\tilde{\sigma}, \tilde{\rho})&=
\beta_g\left(T\tilde{\rho}\beta_g^{-1}(V_\sigma)\right)V_\sigma^*
= \beta_g\left(\tilde{\rho}'\beta_g^{-1}(V_\sigma)T \right)V_\sigma^* \\
&= \beta_g\tilde{\rho}'\beta_g^{-1}(V_\sigma)\beta_g\tilde{\bar{\sigma}}^{\mathrm{op}}(T) V_\sigma^*
= \beta_g\tilde{\rho}'\beta_g^{-1}(V_\sigma)V_\sigma^*\tilde{\sigma}(T) \\
&= \mathcal{E}(\tilde{\sigma}, \tilde{\rho}')\tilde{\sigma}(T).
\end{align*}

(BF2).
Let $h[\tilde{\xi}]=\Ad u \circ \tilde{\xi}^{(h)}$. By Proposition \ref{prop:eta-inter},
$u^*\mathcal{E}(\tilde{\rho},\tilde{\xi})\in (\tilde{\rho}\tilde{\xi},\tilde{\xi}^{(h)}\tilde{\rho})\subset
(\rho\xi,\xi^{(h)}\rho)\otimes \mathbb{C}$.
Thus we have $\tilde{\bar{\sigma}}^{\mathrm{op}}(u^*\mathcal{E}(\tilde{\rho},\tilde{\xi}))
=u^*\mathcal{E}(\tilde{\rho},\tilde{\xi})$, and
\[
V_\sigma \beta_g\left(u^*\mathcal{E}(\tilde{\rho},\tilde{\xi}) \right)V_\sigma^*
=V_\sigma \beta_g\tilde{\bar{\sigma}}^{\mathrm{op}}\left(u^*\mathcal{E}(\tilde{\rho},\tilde{\xi})\right) V_\sigma^*
=\tilde{\sigma}\left(u^*\mathcal{E}(\tilde{\rho},\tilde{\xi})\right).
\]
So
\[
\tilde{\sigma}\left(\mathcal{E}(\tilde{\rho},\tilde{\xi})\right)=
\tilde{\sigma}(u)V_\sigma \beta_g\left(u^*\mathcal{E}(\tilde{\rho},\tilde{\xi})\right) V_\sigma^*
\]
holds.
Also note $\mathcal{E}(\tilde{\sigma},h[\tilde{\xi}])=\beta_g(u)
\mathcal{E}(\tilde{\sigma},\tilde{\xi}^{(h)})\tilde{\sigma}(u^*)$.
Then the first relation of (BF2) can be verified as follows;
\begin{align*}
&\mathcal{E}(\tilde{\sigma},h[\tilde{\xi}])\tilde{\sigma}\left(\mathcal{E}(\tilde{\rho},\tilde{\xi})\right)
=\mathcal{E}(\tilde{\sigma},h[\tilde{\xi}])
\tilde{\sigma}(u)V_\sigma \beta_g \left( u^*\mathcal{E}(\rho,\xi) \right) V_\sigma^* \\
&=\beta_g(u)
\mathcal{E}(\tilde{\sigma},\tilde{\xi}^{(h)})
V_\sigma\beta_g \left( u^*\mathcal{E}(\rho,\xi) \right) V_\sigma^* 
=\beta_g(u)\beta_g\tilde{\xi}^{(h)}\beta_g^{-1}(V_\sigma)
\beta_g\left(u^*\mathcal{E}(\rho,\xi) \right) V_\sigma^*\\
&=\beta_{gh}\tilde{\xi}\beta_{gh}^{-1}(V_\sigma)
\beta_{gh}\tilde{\xi}\beta_h^{-1}(V_\rho)\beta_{g}(V_\rho^*)V_\sigma^* \\
&=\beta_{gh}\tilde{\xi}\beta_{gh}^{-1}\left(V_\sigma\beta_g(V_\rho)\right)V_{\sigma\rho}^* 
=\beta_{gh}\tilde{\xi}\beta_{gh}^{-1}\left(V_{\sigma\rho}\right)V_{\sigma\rho}^*\,\, 
(\mbox{by Proposition } \ref{prop:Vunitary} (\mathrm{ii})) \\
&=\mathcal{E}(\tilde{\sigma}\tilde{\rho},\xi).
\end{align*}

Take $T\otimes 1\in (\tilde{\sigma}, \tilde{\sigma}')$ with  $T\in (\sigma,\sigma')$.
By Proposition \ref{prop:Vunitary}, we have 
$(T\otimes 1)V_\sigma=V_{\sigma'}\beta_g\left(1\otimes j(\hat{T})\right)$ and 
$(T^*\otimes 1)V_{\sigma'}=V_\sigma\beta_g\left(1\otimes j(\hat{T}^*)\right)$.

Hence
\begin{align*}
g[\tilde{\xi}](T\otimes 1)\mathcal{E}(\tilde{\sigma}, \tilde{\xi})&=
\beta_g\tilde{\xi}\beta_g^{-1}\left((T\otimes 1)V_\sigma\right)V_\sigma^*
=
\beta_g\tilde{\xi}\beta_g^{-1}\left(V_{\sigma'}\beta_g(1\otimes j(\hat{T}))\right)V_\sigma^* \\
&=
\beta_g\tilde{\xi}\beta_g^{-1}\left(V_{\sigma'}\right)\beta_g(1\otimes j(\hat{T}))V_{\sigma}^*
=
\beta_g\tilde{\xi}\beta_g^{-1}\left(V_{\sigma'}\right)V_{\sigma'}^*(T\otimes 1) \\
&=\mathcal{E}(\tilde{\sigma}', \tilde{\xi}) (T\otimes 1).
\end{align*}

(BF3).
Let
$v\in (\tilde{\rho}^{(k)},k[\tilde{\rho}])$ and
$w\in (\tilde{\bar{\sigma}}^{(k),\mathrm{op}}, k[\tilde{\bar{\sigma}}^{\mathrm{op}}])$.
Put $u:=\beta_{k}(V_\sigma)\beta_{kgk^{-1}}(w)V_{\sigma^{(k)}}^*$. Then
$k[\tilde{\sigma}]=\Ad u\circ \tilde{\sigma}^{(k)}$ holds, and  we have
\begin{align*}
\mathcal{E}(k[\tilde{\sigma}],k[\tilde{\rho}])&=
\mathcal{E}(\Ad u\circ \tilde{\sigma}^{(k)},\Ad v\circ \tilde{\rho}^{(k)}) \\
&=\beta_{kgk^{-1}}(v)\beta_{kgk^{-1}}\tilde{\rho}^{(k)}\beta_{kgk^{-1}}^{-1}(u)
\mathcal{E}(\tilde{\sigma}^{(k)},\tilde{\rho}^{(k)})
u^*k(\tilde{\sigma})(v^*) \\
&=
\beta_{kgk^{-1}}(v)\beta_{kgk^{-1}}\tilde{\rho}^{(k)}\beta_{kgk^{-1}}^{-1}(u)
\beta_{kgk^{-1}}\tilde{\rho}^{(k)}\beta_{kgk^{-1}}^{-1}(V_{\sigma^{(k)}})V_{\sigma^{(k)}}^*
u^*k(\tilde{\sigma})(v^*) \\
&=
\beta_{kgk^{-1}}(v)
\beta_{kgk^{-1}}\tilde{\rho}^{(k)}\beta_{kgk^{-1}}^{-1}\left(uV_{\sigma^{(k)}}\right)
V_{\sigma^{(k)}}^*
u^*k(\tilde{\sigma})(v^*) \\
&=(*).
\end{align*}
By the definition of $u$, $uV_{\sigma^{(k)}}=\beta_{k}(V_\sigma)\beta_{kgk^{-1}}(w)$.
So
\begin{align*}
(*)&=
\beta_{kgk^{-1}}(v)\beta_{kgk^{-1}}\tilde{\rho}^{(k)}(\beta_{kg^{-1}}(V_\sigma)w)
\beta_{kgk^{-1}}\left(w^*\right)\beta_{k}(V_\sigma^*)
k[\tilde{\sigma}](v^*) \\
&=\beta_{kg}\tilde{\rho}\beta_g^{-1}(V_\sigma)\beta_{kgk^{-1}}(v\tilde{\rho}^{(k)}(w))
\beta_{kgk^{-1}}\left(w^*\right)\beta_{k}(V_\sigma^*)
k[\tilde{\sigma}](v^*) \\
&=
\beta_{kg}\tilde{\rho}\beta_g^{-1}(V_\sigma)
\beta_{kgk^{-1}}(v\tilde{\rho}^{(k)}(w))
\beta_{kgk^{-1}}\left(w^*\right)\beta_{k}(V_\sigma^*\tilde{\sigma}\beta_k^{-1}(v^*) )\\
&=
\beta_{kg}\tilde{\rho}\beta_g^{-1}(V_\sigma)
\beta_{kgk^{-1}}(v\tilde{\rho}^{(k)}(w))
\beta_{kgk^{-1}}\left(w^*\right)\beta_{k}(\beta_g\tilde{\bar{\sigma}}^{\mathrm{op}}\beta_k^{-1}(v^*)V_\sigma^* )\\
&=
\beta_{kg}\tilde{\rho}\beta_g^{-1}(V_\sigma)
\beta_{kgk^{-1}}\left(v\tilde{\rho}^{(k)}(w)\tilde{\bar{\sigma}}^{(k),\mathrm{op}}(v^*)w^*\right)
\beta_{k}(V_\sigma^* ) \\
&=
\beta_{kg}\tilde{\rho}\beta_g^{-1}(V_\sigma)
\beta_{k}(V_\sigma^* ) \,\,\, (\mbox{by Lemma \ref{lem:Intercommute}})\\
&=\beta_k(\mathcal{E}(\tilde{\sigma},\tilde{\rho})).
\end{align*}
\hfill$\Box$

\medskip

At the end of this section, we discuss another type of $G$-braiding. 
If we examine the proof of Theorem \ref{thm:Gbraid1}, we notice the existence of 
another $G$-braiding $\mathcal{E}^{-}(\tilde{\sigma},\tilde{\rho})$.
(To distinguish two $G$-braidings, we denote the unitary $V_\sigma$ and 
the $G$-braiding constructed in Theorem \ref{thm:Gbraid1}
by $V_\sigma^+$ and $\mathcal{E}^+(\tilde{\sigma},\tilde{\rho})$, respectively.)

\begin{lem}\label{lem:Vminus0}
Let $V_\xi^{-}:=(1\otimes j(R_\xi^*))(\id\otimes \bar{\xi}^{\mathrm{op}})\left(\kappa_{{\xi}}U_h^*\right)\in B(\fH_0) $,
$\xi\in \cC_h$.
Then $V_\xi^{-}$ is a unitary with $V_\xi^{-*}=(\bar{R}_\xi^*\otimes 1)(\xi\otimes \id)(\kappa_{\bar{\xi}}U_{h^{-1}}^*) $, 
and $V_\xi^{-} \bar{\xi}^{\mathrm{op}}\gamma_L(\beta_h(a))=\xi\gamma_L(a)V_\xi^{-}$ for $a\in \cB_0$.
\end{lem}
\textbf{Proof.} This can be shown by Lemma \ref{lem:LR-relation}. Note that $U_h\in \gamma_L(\cB)$.
We remark that $V_\xi^{-}$ for $d(\xi)=\infty$ should be understood as 
$\sum_{k}(w_k\otimes 1)V_{\pi_k}^-(1\otimes j(\bar{w}_k))$ via  
irreducible decomposition $\xi(x)=\sum_k w_k\pi_k(x)w_k^*$, $\bar{\xi}(x)=\sum_k \bar{w}_k\bar{\pi}_k(x)\bar{w}_k^*$.
\hfill$\Box$

\begin{lem}\label{lem:Vminus1}
Define
\[
 \tilde{V}_{\bar{\rho}}^{-}:=(1\otimes j(R_\rho^*))\tilde{\bar{\rho}}^{\mathrm{op}}(\hlambda_\rho),\,\,\, 
\mathcal{E}^-(\tilde{\sigma},\tilde{\rho}):=\tilde{V}_{\tilde{\rho}}^-\tilde{\sigma}(\tilde{V}_{\tilde{\rho}}^{-})^*,\,\,
\tilde{\sigma} \in Z_{\cC_0}(\cC),  \,\, \tilde{\rho} \in Z_{\cC_0}(\cC)_h.
\]
Then $\tilde{V}_{\tilde{\rho}}^-$ is a unitary
in $
\left(\tilde{\bar{\rho}}^{\mathrm{op}}\beta_h, \tilde{\rho}\right)$.
Hence $\mathcal{E}^-(\tilde{\sigma},\tilde{\rho})$  is also a unitary and 
$\mathcal{E}^-(\tilde{\sigma},\tilde{\rho}) \in (\tilde{\sigma}\tilde{\rho}, \tilde{\rho} h^{-1}[\tilde{\sigma}])$.
\end{lem}
\textbf{Proof.}
We will show $ \cE_\rho V_{\rho}^{-}
\cE_{\bar{\rho}}^{\mathrm{op}*}=\gamma_L(\tilde{V}_{\tilde{\rho}}^-)$.
Then we obtain conclusion by Lemma \ref{lem:Vminus0}.
(It seems difficult to see the unitarity of $\tilde{V}_{\tilde{\rho}}^-$ directly.)
At first we compute  $\cE_\rho(1\otimes j(R_\rho^*))(\id\otimes \bar{\rho}^{\mathrm{op}})(\kappa_\rho)$ as follows.
Take any $v\in \fH$. Then
\begin{align*}
&\left( \cE_\rho(1\otimes j(R_\rho^*))(\id\otimes \bar{\rho}^{\mathrm{op}})(\kappa_\rho)v\right)(\xi) 
=
\sum_{\eta\in \cC_0}(\cE_\rho(\xi)\otimes j(R_\rho^*))(\id\otimes \bar{\rho}^{\mathrm{op}})(T_{\rho,\xi}^{\eta})v(\eta) \\
&=
\sum_{\eta\in \cC_0}(1 \otimes j(R_\rho^*))(\id\otimes \bar{\rho}^{\mathrm{op}})((\cE_\rho(\xi)\otimes 1)T_{\rho,\xi}^{\eta})v(\eta) \\
&=
\sum_{\eta\in \cC_0}(1 \otimes j(R_\rho^*))(\id\otimes \bar{\rho}^{\mathrm{op}})
((1\otimes j(\cE_\rho(\xi))^*)T_{\xi, \rho}^{\eta})v(\eta) \\
&=
(1 \otimes j(R_\rho^*\bar{\rho}(\cE_\rho(\xi))^*))
(\id\otimes \bar{\rho}^{\mathrm{op}})(\gamma_L(\lambda_\rho)v)(\xi) \\
&=
(1 \otimes j(\xi(R_\rho^*)\cE_{\bar{\rho}}(\xi)))
(\id\otimes \bar{\rho}^{\mathrm{op}})(\gamma_L(\lambda_\rho)v)(\xi) \\
&=\gamma_L(1\otimes j(R_\rho)^*)
\cE_{\bar{\rho}}(\id\otimes \bar{\rho}^{\mathrm{op}})(\gamma_L(\lambda_\rho)v)(\xi).
\end{align*}
Hence we have
\begin{align*}
\cE_\rho V_{\rho}^{-}\cE_{\bar{\rho}}^{\mathrm{op}*}&=
 \cE_\rho(1\otimes j(R_\rho^*))(\id\otimes \bar{\rho}^{\mathrm{op}})(\kappa_\rho U_h^*)
\cE_{\bar{\rho}}^{\mathrm{op}*} \\
&=
\gamma_L(1\otimes j(R_\rho)^*)\cE_{\bar{\rho}}(\id\otimes \bar{\rho}^{\mathrm{op}})(\gamma_L(\lambda_\rho)
U_h^*)
\cE_{\bar{\rho}}^{\mathrm{op}*} \\
&=\gamma_L(1\otimes j(R_\rho)^*\tilde{\bar{\rho}}^{\mathrm{op}}(\hlambda_\rho)) 
=\gamma_L(\tilde{V}_{\tilde{\rho}}^-).
\end{align*}
\hfill$\Box$

\begin{thm}\label{thm:Gbraid+-}
 We have
$\mathcal{E}^-(\tilde{\sigma},\tilde{\rho})=\mathcal{E}^+(\tilde{\rho},h^{-1}[\tilde{\sigma}])^*$.
\end{thm}
\textbf{Proof.} We can see $\tilde{\rho}=\Ad\left(\tilde{V}_{\tilde{\rho}}^{-}\beta_h^{-1}(V_\rho^{+*})\right)\circ 
h^{-1}[\tilde{\rho}]$.
By (BF2) and (BF3), 
\begin{align*}
 \mathcal{E}^+(\tilde{\rho},h^{-1}[\tilde{\sigma}])
&=
\tilde{\sigma}(\tilde{V}_{\tilde{\rho}}^-\beta_h^{-1}(V_{{\rho}}^{+*}))
 \mathcal{E}^+\left(h^{-1}(\tilde{\rho}),h^{-1}(\tilde{\sigma})\right)
\beta_h^{-1}(V_{{\rho}}^+)V_{\tilde{\rho}}^{-*} \\
&=
\tilde{\sigma}(V_{\tilde{\rho}}^-\beta_h^{-1}(V_\rho^{+*}))
 \beta_h^{-1}\left(\mathcal{E}^+\left(\tilde{\rho},\tilde{\sigma}\right)\right)
\beta_h^{-1}(V_\rho^+)\tilde{V}_{\tilde{\rho}}^{-*} \\
&=
\tilde{\sigma}(\tilde{V}_{\tilde{\rho}}^-\beta_h^{-1}(V_\rho^{+*}))
 \beta_h^{-1}\left(\beta_h\tilde{\sigma}\beta_h^{-1}(V_\rho^+)V_\rho^{+*}\right)
\beta_h^{-1}(V_\rho^+)\tilde{V}_{\tilde{\rho}}^{-*} \\
&=
\tilde{\sigma}(\tilde{V}_{\tilde{\rho}}^-)\tilde{V}_{\tilde{\rho}}^{-*} 
=\mathcal{E}^-(\tilde{\sigma},\tilde{\rho})^*.
\end{align*}
\hfill$\Box$

\section{$G$-equivariant categories and center  for full categories}\label{sec:full}
We keep the notation from the previous section.
\begin{thm}\label{thm:QsysFull}
An inclusion  $\cA\subset \cB_0\rtimes_\beta G$ is isomorphic to 
the Longo-Rehren inclusion $\cA\subset \cB$ for $\cC$.
\end{thm}
\textbf{Proof.} 
Define  $\bar{U}$ be a unitary in $B(\fH, \ell^2(G,\fH_0))$ by 
$\bar{U}(\oplus_g v_g)=\oplus_g U_gv_g$. Then we can easily see 
$\bar{U}\gamma_L(a)\bar{U}^*=\bigoplus_g \gamma_L(\beta_g(a))$, 
$a\in \cB_0$, 
and $\bar{U} U_g\bar{U}^*=1\otimes \lambda^G_g$, where $\lambda^G_g$ is the  right regular representation of $G$ on 
$\ell^2(G)$. This shows the statement.
\hfill$\Box$

\medskip

Since $U_g\in \gamma_L(\cB)$, we 
can take $u_g^G\in \cB$ with $\gamma_L(u_g^G)=U_g$. By Theorem \ref{thm:QsysFull}, 
we can regard $u_g^G$ as an implementing unitary for $\cB=\cB_0\rtimes_\beta G$.

\begin{df}\label{df:fixed}
$(1)$ We say $\theta\in Z_{\cC_0}(\cC)$ is $G$-equivariant if 
there exists a family of unitaries $\{c^\theta_g\}\subset \cB_0$ such that 
\[
 \beta_g\theta\beta_g^{-1}=\Ad c^\theta_g\circ \theta,\,\,\,\beta_g(c^\theta_h)c^\theta_g=c^\theta_{gh}. 
\]
$(2)$ Let $(\theta,\{c_g^\theta\}_{g\in G})$ be a $G$-equivariant element. The 
extension $\hat{\eta}(\theta, c^\theta)\in \End(\cB_0\rtimes_\beta G)$ of $\theta$ is defined by
by $\hat{\eta}(\theta,c^\theta)(u_g^G):=
c_g^{\theta*}u_g^G$. \\
$(3)$ 
Define  $Z_{\cC_0}(\cC)^G\subset \End(\cB_0\rtimes_\beta G)$ by
\[
Z_{\cC_0}(\cC)^G:=  \{\mu\in \End_0(\cB_0\rtimes_\beta G)\mid \mu \sim \hat{\eta}(\theta,c^\theta) \mbox{ for some 
$G$-equivariant element }
(\theta,c_g^\theta)\}.
\]
\end{df}

We can show the existence of the extension $\hat{\eta}(\theta, c^\theta)$ in a similar
way as that of $\eta$-extension.

Let $\iota_G$ is the inclusion map $\cB_0\subset \cB_0\rtimes_\beta G$.
We can identify  $Z_{\cC_0}(\cC)^G$ with
 \begin{align*}
&
\{\mu\in \End_0(\cB_0\rtimes_\beta G)\mid \mu\iota_G \sim \iota_G\theta \mbox{ for some }
\theta \in Z_{\cC_0}(\cC)\} \\
=&
\{\mu\in \End_0(\cB_0 \rtimes_\beta G)\mid \mu \prec \iota_G\theta \bar{\iota}_G\mbox{ for some }
\theta \in Z_{\cC_0}(\cC)\}\,\,  (\mbox{when } |\Irr(\cC)|<\infty ).
 \end{align*}
To see the last equality, one should notice the irreducibility of  $\iota_G\theta$  for $\theta\in \mathrm{Irr}(Z_{\cC_0}(\cC))$,
and $\iota_G\beta_g\sim \iota_G$. Then it is shown similarly as in \cite[Theorem 3.2]{Kw-Dcenter}.

In the following, we often 
denote $\hat{\eta}(\tilde{\sigma},c^{\tilde{\sigma}})$ by $\hat{\sigma}$ simply.
In a similar way as in \cite[Theorem 4.6]{Iz-LR}, we can see the following. (Also see Proposition \ref{prop:eta-inter}.)
\begin{lem}\label{lem:hateta-inter}
$(1)$
\[
 \left(\hat{\eta}(\theta, c^\theta),\hat{\eta}(\theta',c^{\theta'})\right)=
\{T\in (\theta,\theta')\mid \beta_g(T)c_g^{\theta}=c_g^{\theta'}T, g\in G\}.
\]
$(2)$ For $(\theta,c_g^\theta)\in Z_{\cC_0}(\cC)^G$, define $c_g^{\bar{\theta}}$ by 
$c_g^{\bar{\theta}}:=R_\theta^* \bar{\theta}(c_g^{{\theta}*})\bar{\theta}\beta_g(\bar{R}_\theta)$.
Then $c_g^{\bar{\theta}}$ is a unitary with $(\bar{\theta},c_g^{\bar{\theta}})\in Z_{\cC_0}(\cC)^G$,
and $\overline{\hat{\eta}(\theta,c^\theta)}=\hat{\eta}(\bar{\theta},c^{\bar{\theta}})$ holds.
\end{lem}

\begin{lem}\label{lem:fixedHalf}
Let $\tilde{\sigma}\in Z_{\cC_0}(\cC)$ so that $\beta_g\tilde{\sigma}\beta_g^{-1}\sim \tilde{\sigma}$ for 
some $g\in G$,
and fix a unitary $c_g^{\tilde{\sigma}}\in \cB_0$ such that $\beta_g\tilde{\sigma}\beta_g^{-1}=\Ad c_g^{\tilde{\sigma}} 
\tilde{\sigma}$.
Let $\{\mathcal{E}_\sigma(\pi)\}_{\pi\in \cC_g}$ be a family of unitaries 
for $c_g^{\tilde{\sigma}}$ given in Lemma \ref{lem:actinInter}.
Then the equation $\beta_g(c_h^{\tilde{\sigma}})c_g^{\tilde{\sigma}}=c_{gh}^{\tilde{\sigma}}$ 
is equivalent to the half braiding property of 
$\{\cE_\sigma(\pi)\}_{\pi\in \cC}$. 
\end{lem}
\textbf{Proof.}
In this proof, we denote the restriction of $\gamma_L(a)$ on $\fH_g$ by $\gamma_{L}^g(a)$.

Let $\cE_\sigma^g \in \End(\fH_{g,\cB_0})$ by $(\cE_\sigma^g v)(\pi)=(\cE_\sigma(\pi)\otimes 1)v(\pi)$.
By the construction in Lemma \ref{lem:actinInter}, we have
\[
\gamma_L^0(c_g^{\tilde{\sigma}})= U_g\cE_\sigma^g(\sigma\otimes \id)(U_g^*)\cE_\sigma^*.
\]
\begin{align*}
 \beta_g(c_h^{\tilde{\sigma}})c_g^{\tilde{\sigma}}=c_{gh}^{\tilde{\sigma}}
&\Leftrightarrow 
 \gamma_{l}^0\left(\beta_g(c_h^{\tilde{\sigma}})c_g^{\tilde{\sigma}}\right)=\gamma_L^0(c_{gh}^{\tilde{\sigma}}) \\
&\Leftrightarrow  
U_g \gamma_{l}^g\left(c_h^{\tilde{\sigma}}\right)U_g^*
U_g\cE_\sigma^g(\sigma\otimes \id)(U_g^*)\cE_\sigma^*.
=U_{gh}\cE_\sigma^{gh}(\sigma\otimes \id)(U_{gh}^*)\cE_\sigma^* \\
&\Leftrightarrow  
\gamma_{l}^g\left(c_h^{\tilde{\sigma}}\right)
=U_{h}\cE_\sigma^{gh}(\sigma\otimes \id)(U_{h}^*)\cE_\sigma^{g*}
\end{align*}

Take any $\pi\in \Irr(\cC_g)$.
Since $\kappa_\pi^*\kappa_\pi=d(\pi)$,
the last equation is equivalent to 
\[
\kappa_\pi\gamma_{l}^g\left(c_h^{\tilde{\sigma}}\right)
=\kappa_\pi U_{h}\cE_\sigma^{gh}(\sigma\otimes \id)(U_{h}^*)\cE_\sigma^{g*}.
\]
The left hand side is 
\[
\kappa_\pi \gamma_{l}^g\left(c_h^{\tilde{\sigma}}\right)=
(\pi\otimes \pi^\mathrm{op})\gamma_{l}^0\left(c_h^{\tilde{\sigma}}\right)\kappa_\pi=
 (\pi\otimes \pi^\mathrm{op})\left(
U_h\cE_\sigma^h(\sigma\otimes \id)(U_h^*)\cE_\sigma^*\right)\kappa_\pi.
\]
By Lemma \ref{lem:actinInter}, we have
\[
 (\pi\otimes \pi^\mathrm{op})(\cE_\sigma^*)
\kappa_\pi\cE_\sigma^{g}=
(\cE_\sigma(\pi)\otimes 1)(\sigma\otimes \id)(\kappa_\pi).
\]

Thus we have
\begin{align*}
&\kappa_\pi\gamma_{l}^g\left(c_h^{\tilde{\sigma}}\right)
=\kappa_\pi U_{h}\cE_\sigma^{gh}(\sigma\otimes \id)(U_{h}^*)\cE_\sigma^{g*} \\
&\Leftrightarrow 
 (\pi\otimes \pi^\mathrm{op})\left(
\cE_\sigma^h\right)(\pi\sigma\otimes \pi^{\mathrm{op}}) (U_h^*)
(\pi\otimes \pi^\mathrm{op})(\cE_\sigma^*)
\kappa_\pi\cE_\sigma^{g}
 =\kappa_\pi \cE_\sigma^{gh}(\sigma\otimes \id)(U_{h}^*) \\
&\Leftrightarrow 
 (\pi\otimes \pi^\mathrm{op})\left(
\cE_\sigma^h\right)(\pi\sigma\otimes \pi^{\mathrm{op}}) (U_h^*)
(\cE_\sigma(\pi)\otimes 1)(\sigma\otimes \id)(\kappa_\pi)
 =\kappa_\pi \cE_\sigma^{gh}(\sigma\otimes \id)(U_{h}^*) \\
&\Leftrightarrow 
 (\pi\otimes \pi^\mathrm{op})\left(
\cE_\sigma^h\right)(\cE_\sigma(\pi)\otimes 1)
(\sigma\pi\otimes \pi^{\mathrm{op}}) (U_h^*)
(\sigma\otimes \id)(\kappa_\pi)
 =\kappa_\pi \cE_\sigma^{gh}(\sigma\otimes \id)(U_{h}^*) \\
&
\Leftrightarrow 
 (\pi\otimes \pi^\mathrm{op})\left(
\cE_\sigma^h\right)(\cE_\sigma(\pi)\otimes 1)
(\sigma\otimes \id)(\kappa_\pi)
(\sigma\otimes \id ) (U_h^*)
 =\kappa_\pi \cE_\sigma^{gh}(\sigma\otimes \id)(U_{h}^*) \\
&\Leftrightarrow 
 (\pi\otimes \pi^\mathrm{op})\left(
\cE_\sigma^h\right)(\cE_\sigma(\pi)\otimes 1)
(\sigma\otimes \id)(\kappa_\pi)
 =\kappa_\pi \cE_\sigma^{gh}.
\end{align*}
The last equation is equivalent to the half braiding property of $\cE_\sigma(\cdot)$.
\hfill$\Box$

\begin{prop}\label{prop:ext-eta-full}
 For $(\tilde{\sigma}, c_g^{\tilde{\sigma}})\in Z_{\cC_0}(\cC)^G$, the extension $\hat{\sigma}$ 
of $\tilde{\sigma}$ to $\cB_0 \rtimes_\beta G$ is equal to the $\eta$-extension of $(\sigma\otimes \id,\cE_\sigma)$ 
for $\cB$. Here $\cE_\sigma$ is a half braiding  corresponding to $c_g^{\tilde{\sigma}}$ in 
Lemma \ref{lem:fixedHalf}.
\end{prop}
\textbf{Proof.}
By the definition of $\hat{\sigma}$, we have
\begin{align*}
 \hat{\sigma}(\lambda_\pi)&=\hat{\sigma}(\hlambda_\pi u^G_g) = 
\iota_G(\tilde{\sigma}(\hlambda_\pi)c_g^{\tilde{\sigma}*})u_g^G.
\end{align*}

Here
\begin{align*}
\gamma_L(\tilde{\sigma}(\hlambda_\pi)c_g^{\tilde{\sigma}*})&=
\cE_\sigma \sigma\gamma_L(\hlambda_\pi)\cE_\sigma^* \cE_\sigma\sigma(U_g)\cE_\sigma^*
U_g^* \\
&=
\cE_\sigma \sigma\left(\gamma_L(\hlambda_\pi)U_g\right)\cE_\sigma^*U_g^*
=
\cE_\sigma \sigma\left(\gamma_L(\lambda_\pi) \right)\cE_\sigma^*U_g^* \\
&=\gamma_L((\cE_\sigma(\pi)^*\otimes 1)\lambda_\pi  )U_g^* =\gamma_L((\cE_\sigma(\pi)^*\otimes 1)\hlambda_\pi  )
\end{align*}

Thus
$\hat{\sigma}(\lambda_\pi)
=
\left(\cE_\sigma(\pi)^*\otimes 1\right)\lambda_\pi$.
\hfill$\Box$

\bigskip

 Lemma \ref{lem:fixedHalf} and Proposition \ref{prop:ext-eta-full} imply the following
\cite[Theorem 3.5]{GNN-center}.

\begin{thm}\label{thm:full-fix-indentify}
We can identify 
 $Z_{\cC_0}(\cC)^G$ and $Z(\cC)$.
\end{thm}

\begin{lem}\label{lem:fixbraid0}
 Let $\hat{\cE}(\tilde{\sigma},\hat{\rho}):=c_g^{\tilde{\rho}*}\cE(\tilde{\sigma},\tilde{\rho})$ 
for $\tilde{\sigma}\in Z_{\cC_0}(\cC)_g$, $(\tilde{\rho},c_g^{\tilde{\rho}})\in Z_{\cC_0}(\cC)^G$, and 
extend it canonically for every $\tilde{\sigma}\in Z_{\cC_0}(\cC)$.
Then $\hat{\cE}(\tilde{\sigma},\hat{\rho})=\cE_{\rho}(\sigma)^*\otimes 1$.
\end{lem}
\textbf{Proof.} It suffices to show the statement for $\tilde{\sigma}\in Z_{\cC_0}(\cC)_g$.
By $\hat{\rho}(u_g^G)=c_g^{\tilde{\rho}*}u_g^G$,
\begin{align*}
\hat{\cE}(\tilde{\sigma},\hat{\rho})= c_g^{\tilde{\rho}*}\cE(\tilde{\sigma},\tilde{\rho})&=
c_g^{\tilde{\rho}*}\beta_g\tilde{\rho}\beta_g^{-1}(V_\sigma)V_\sigma^* 
=
\tilde{\rho}(V_\sigma)c_g^{\tilde{\rho}*}V_\sigma^* 
=
\hat{\rho}(V_\sigma u_g^G) u_g^{G*}V_\sigma^* 
\end{align*}
holds. Therefore
\[
\hat{\rho}(V_\sigma u_g^G)u _g^{G*}V_\sigma^* =
\hat{\rho}((1\otimes j(\bar{R}_\sigma^*))\lambda_\sigma)u_g^{G*} V_\sigma^*
=(\cE_\rho(\sigma)^*\otimes 1) V_\sigma V_\sigma^*
=\cE_\rho(\sigma)^*\otimes 1
\]
by Proposition \ref{prop:ext-eta-full}. 
\hfill$\Box$

\bigskip

Since $\cE_\rho(\sigma)^*\otimes 1$ gives a braiding on $Z(\cC)$ (see \cite[Proposition 6.5]{Iz-LR}), 
$\hat{\cE}$ also gives a braiding on $Z_{\cC_0}(\cC)^G$.
In the following, we will explain this fact directly without referring to $Z(\cC)$.

\begin{lem}\label{lem:fixedbraid1}
For $\tilde{\sigma}\in {Z}_{\cC_0}(\cC)$ and $\hat{\rho}\in Z_{\cC_0}(\cC)^G$,
$\hat{\cE}(\tilde{\sigma},\hat{\rho})\in (\tilde{\sigma}\tilde{\rho},\tilde{\rho}\tilde{\sigma})$ and
 \begin{align*}
 T\hat{\cE}(\tilde{\sigma},\hat{\rho})&= \hat{\cE}(\tilde{\sigma},\hat{\rho}')\tilde{\sigma}(T), 
\,\,\, T\in (\hat{\rho},\hat{\rho}'),
\\
\hat{\cE}(\tilde{\sigma},\hat{\rho}_1\hat{\rho}_2) &= 
\hat{\rho}_1\left(\hat{\cE}(\tilde{\sigma},\hat{\rho}_2)\right)  \hat{\cE}(\tilde{\sigma},\hat{\rho}_1), \\
\hat{\rho}(T)\hat{\cE}(\tilde{\sigma},\hat{\rho})&=
\hat{\cE}(\tilde{\sigma}',\hat{\rho})T,\,\, T\in (\tilde{\sigma},\tilde{\sigma}'), \\
\hat{\cE}(\tilde{\sigma}_1\tilde{\sigma}_2,\hat{\rho}) &=
\hat{\cE}(\tilde{\sigma}_1,\hat{\rho}) \tilde{\sigma}_1\left(\hat{\cE}(\tilde{\sigma}_2,\hat{\rho})\right). 
 \end{align*}
\end{lem}
\textbf{Proof.} By Lemma \ref{lem:hateta-inter}, 
$T\in (\hat{\rho},\hat{\rho}')$ if and only if $T\in (\tilde{\rho},\tilde{\rho}')$ and 
$\beta_g(T)c_g^{\tilde{\rho}}=c_g^{\tilde{\rho}'}T$ for all $g\in G$.
Then it is routine to show  statements  by (BF0),  (BF1) and (BF2) of Theorem \ref{thm:Gbraid1}.
\hfill$\Box$

\begin{lem}\label{lem:fixedbraid2}
For $\hat{\sigma},\hat{\rho}\in Z_{\cC_0}(\cC)^G$, define
$\hat{\cE}(\hat{\sigma},\hat{\rho}):=\hat{\cE}(\tilde{\sigma},\hat{\rho})$.
Then we have
$ \hat{\cE}(\hat{\sigma},\hat{\rho})\in (\hat{\sigma}\hat{\rho},\hat{\rho}\hat{\sigma})$.
\end{lem}
\textbf{Proof.}
We already know $ \hat{\cE}(\hat{\sigma},\hat{\rho})\in (\tilde{\sigma}\tilde{\rho},\tilde{\rho}\tilde{\sigma})$
in Lemma \ref{lem:fixedbraid1}.
By Lemma \ref{lem:hateta-inter}, we only have to show
\[
{\beta_g\left(\hat{\cE}(\hat{\sigma},\hat{\rho})\right)
c_g^{\tilde{\sigma}}\tilde{\sigma}(c_g^{\tilde{\rho}}) } =
c_g^{\tilde{\rho}}\tilde{\rho}(c_g^{\tilde{\sigma}})\hat{\cE}(\hat{\sigma},\hat{\rho}).
\]

There exists $L\subset G$ such that  
$\sigma\sim \bigoplus_{l\in L}\sigma_l$, $\sigma_l\in \cC_l$.
Note that this decomposition is unique up to unitary equivalence. 
Since $\beta_g \tilde{\sigma}\beta_g^{-1}\sim \tilde{\sigma}$ for all $g\in G$ 
and $\beta_g\tilde{\sigma}_l\beta_g^{-1}\in Z_{\cC_0}(\cC)_{glg^{-1}}$, we have
$\beta_g\tilde{\sigma}_l\beta_g^{-1}\sim \tilde{\sigma}_{glg^{-1}}$ and hence
$gLg^{-1}=L$ for all $g\in G$.

Take   a family of orthogonal isometries $\{w_l\}_{l\in L}\subset \cM$ with support 1 such that 
$\sigma(x)=\sum_{l\in L}w_l\sigma_l(x)w_l^*$, and set 
\[
c_{g,l}:=\beta_g(w_l^*\otimes 1)c_{g}^{\tilde{\sigma}}
(w_{glg^{-1}}\otimes 1)
\in (\tilde{\sigma}_{glg^{-1}}, \beta_g\tilde{\sigma}_l\beta_g^{-1}).
\]
(In the rest of this proof, we simply denote $w_l\otimes 1$ by $w_l$.)
Since $\beta_g\tilde{\sigma}_l\beta_g^{-1}\in Z_{\cC_0}(\cC)_{glg^{-1}}$, we have
\[
 \beta_g(w_l^*)c_{g}^{\tilde{\sigma}}(w_{gl'g^{-1}})
\in (\tilde{\sigma}_{gl'g^{-1}}, \beta_g\tilde{\sigma}_l\beta_g^{-1})=0
\]
for $l\ne l'$. 
This implies the unitarity of $c_{g,l}$, and 
\[
 \beta_g\tilde{\sigma}_l\beta_g^{-1}=\Ad c_{g,l}\circ \tilde{\sigma}_{glg^{-1}},\,\,
c_g^{\tilde{\sigma}}=\sum_{l\in L}\beta_g(w_l)c_{g,l}w_{glg^{-1}}^*.
\]

On one hand, we have
\begin{align*}
\lefteqn{\beta_g\left(\hat{\cE}(\tilde{\sigma},\hat{\rho})\right)
c_g^{\tilde{\sigma}}\tilde{\sigma}(c_g^{\tilde{\rho}}) } \\
&=\left(\sum_{l}\beta_g\tilde{\rho}(w_l)
 \beta_g\left(\hat{\cE}(\tilde{\sigma}_l,\hat{\rho})\right)\beta_g(w_l^*)\right)
\left(\sum_{l}\beta_g(w_l)
c_{g,l}w_{glg^{-1}}^* \right)\left(\sum_{l}w_l\tilde{\sigma}_l(c_g^{\tilde{\rho}})w_{l}^*\right) \\
&=
\sum_{l}\beta_g\tilde{\rho}(w_l)
 \beta_g\left(\hat{\cE}(\tilde{\sigma}_l,\hat{\rho})\right)
c_{g,l}\tilde{\sigma}_{glg^{-1}}(c_g^{\tilde{\rho}})w_{glg^{-1}}^*  .
\end{align*}

On the other hand, we have
\begin{align*}
c_g^{\tilde{\rho}}\tilde{\rho}(c_g^{\tilde{\sigma}})\hat{\cE}(\tilde{\sigma},\hat{\rho})
&=
c_g^{\tilde{\rho}}\tilde{\rho}\left(\sum_l \beta_g(w_l)c_{g,l}w_{glg^{-1}}^*\right)
\left(\sum_l\tilde{\rho}(w_l)
\hat{\cE}(\tilde{\sigma}_l,\hat{\rho})w_l^*\right)\\
&
=
\sum_l
c_g^{\tilde{\rho}}\tilde{\rho}\left( \beta_g(w_l)c_{g,l}\right)
\hat{\cE}(\tilde{\sigma}_{glg^{-1}},\hat{\rho})w_{glg^{-1}}^*\\
&
=
\sum_l
\beta_g\tilde{\rho}\left( w_l\right)c_g^{\tilde{\rho}}\tilde{\rho}(c_{g,l})
\hat{\cE}(\tilde{\sigma}_{glg^{-1}},\hat{\rho})w_{glg^{-1}}^*.
\end{align*}

Hence we only have to show
\[
 \beta_g\left(\hat{\cE}(\tilde{\sigma}_l,\hat{\rho})\right)
c_{g,l}\tilde{\sigma}_{glg^{-1}}(c_g^{\tilde{\rho}})=c_g^{\tilde{\rho}}\tilde{\rho}(c_{g,l})
\hat{\cE}(\tilde{\sigma}_{glg^{-1}},\hat{\rho})
\]
for every $l\in L$.

By (BF3) of Theorem  \ref{thm:Gbraid1}, 
\begin{align*}
&\beta_g\left(\hat{\cE}(\tilde{\sigma}_{l},\hat{\rho})\right) c_{g,l}
\tilde{\sigma}_{glg^{-1}}(c_g^{\tilde{\rho}})
=
\beta_g(c_{l}^{\tilde{\rho}*})
\beta_g\left({\cE}(\tilde{\sigma}_{l},\tilde{\rho})\right) c_{g,l}
\tilde{\sigma}_{glg^{-1}}(c_g^{\tilde{\rho}}) \\
&=
\beta_g(c_{l}^{\tilde{\rho}*})
{\cE}\left(g(\tilde{\sigma}_{l}),g(\tilde{\rho})\right) c_{g,l}
\tilde{\sigma}_{glg^{-1}}(c_g^{\tilde{\rho}}) 
=\beta_g(c_{l}^{\tilde{\rho}*})
\cE(\Ad c_{g,l}\circ\tilde{\sigma}_{glg^{-1}},\Ad c_g^{\tilde{\rho}}\circ \tilde{\rho}) 
c_{g,l}\tilde{\sigma}_{glg^{-1}}(c_g^{\tilde{\rho}}) \\
&= 
\beta_g(c_{l}^{\tilde{\rho}*})
\beta_{glg^{-1}}(c_g^{\tilde{\rho}})\beta_{glg^{-1}}\tilde{\rho}\beta_{glg^{-1}}(c_{g,l})
\cE(\tilde{\sigma}_{glg^{-1}},\tilde{\rho}) \\
&= 
\beta_g(c_{l}^{\tilde{\rho}*})
\beta_{glg^{-1}}(c_g^{\tilde{\rho}})\beta_{glg^{-1}}\tilde{\rho}\beta_{glg^{-1}}(c_{g,l})
c_{glg^{-1}}^{\tilde{\rho}}
\hat{\cE}(\tilde{\sigma}_{glg^{-1}},\hat{\rho}) 
\\
&= 
\beta_g(c_{l}^{\tilde{\rho}*})
\beta_{glg^{-1}}(c_g^{\tilde{\rho}})
c_{glg^{-1}}^{\tilde{\rho}}
\tilde{\rho}(c_{g,l})
\hat{\cE}(\tilde{\sigma}_{glg^{-1}},\hat{\rho}) 
 = c_g^{\tilde{\rho}}
\tilde{\rho}(c_{g,l})
\hat{\cE}(\tilde{\sigma}_{glg^{-1}},\hat{\rho}).
\end{align*}
The last equation follows from the 1-cocycle property of $c^{\tilde{\rho}}$.
\hfill$\Box$

\bigskip

It is obvious that
Lemma \ref{lem:fixedbraid1} and Lemma \ref{lem:fixedbraid2} implies that $\hat{\cE}$ is a braiding
on $Z_{\cC_0}(\cC)^G$.

\bigskip

We present an example of elements in $Z_{\cC_0}(\cC)^G$. Let $\tilde{\sigma}\in Z_{\cC_0}(\cC)$, and 
 $z_g\in \cB_0$ be a unitary such that $\beta_g\tilde{\sigma} \beta_g^{-1}=\Ad z_g\tilde{\sigma}^{(g)} $. 
Fix a family of orthogonal isometries $\{w_g\}_{g\in G}\subset \cM$ with support 1, and 
set 
\[
 \sigma^G(x):=\sum_{g\in G} w_g\sigma^{(g)}(x)w_g^*,\,\, 
c_g^{\tilde{\sigma},G}:=\sum_{h\in G}\beta_g((w_h\otimes 1)z_h^*)z_{gh}(w_{gh}^*\otimes 1).
\]
Then we can easily see $(\tilde{\sigma}^G, \{c_g^{\tilde{\sigma},G}\})\subset Z_{\cC_0}(\cC)^G$.

We further investigate the above example in two cases.

\noindent
\textbf{Case 1.} $\tilde{\sigma}$ is irreducible, and 
$(\tilde{\sigma}, c_g^{\tilde{\sigma}})\in Z_{\cC_0}(\cC)^G$. (Thus we may assume $\sigma=\sigma^{(g)}$ and 
$z_g=c_g^{\tilde{\sigma}}$.)
Put $u_g:=\sum_{h}w_hw_{hg}^*\otimes 1$. We can see $(\hat{\sigma}^G,\hat{\sigma}^G)
=\{\sum_ga_gu_g\mid a_g\in \mathbb{C}\}$.

Let $G$ be a finite group.
Thus $\hat{\sigma}^G$ can be decomposed according to the irreducible decomposition of the regular representation of $G$.
Indeed, for arbitrary $\chi\in \mathrm{Irr}(G)$, 
fix a family of orthogonal isometries $\{v^\chi_i\}_{1\leq i\leq d\chi}\subset \cM$ with support 1, and 
define 
\[
 \sigma^\chi(x):=\sum_{i=1}^{d(\chi)} v_i^\chi \sigma(x)v_i^{\chi*}, \,\,\,
c_g^{\tilde{\sigma}, \chi}:=\sum_{i,j=1}^{d(\chi)}\beta_g(v_i^{\chi}\otimes 1)\chi(g)_{i,j}c_g^{\tilde{\sigma}}(v_j^{\chi*}\otimes 1).
\]
Then $(\tilde{\sigma}^\chi, c_g^{\tilde{\sigma},\chi})\in Z_{\cC_0}(\cC)^G$, and the irreducible decomposition of 
$\hat{\sigma}^G$ is given by 
$\hat{\sigma}^G= \bigoplus_{\chi\in \mathrm{Irr}(G)}d\chi \hat{\sigma}^\chi$.  
As a special case,  
we can show that  $Z(\cC)$  contains 
$\mathrm{Irr}(G)$ by applying this construction to $\id_{\cB_0}$. 

\noindent
\textbf{Case 2.} $\tilde{\sigma}$ is irreducible, and $\tilde{\sigma}^{(g)}\not\sim \tilde{\sigma}^{(h)}$ for 
$g\ne h$. Note that the choice of $\{z_h\}_{h\in G}$ is not unique.  
Fix  $\{b_h\}_{h\in G}\subset \mathbb{T}$, 
and define $d_g^{\tilde{\sigma}^G}$ by 
\[
d_g^{\tilde{\sigma},G}:=\sum_{h\in G}\beta_g((w_h\otimes 1)\overline{b_h}z_h^*){b_{gh}}z_{gh}(w_{gh}^*\otimes 1).
\]
Let $u=\sum_{h}{b_h}w_hw_h^*$. 
Then we can see
\[
 \left( \hat{\eta}(\tilde{\sigma}^G,d^{\tilde{\sigma},G}) , \hat{\eta}(\tilde{\sigma}^G,c^{\tilde{\sigma},G}) 
\right)=\{c(u\otimes 1)\mid c\in \mathbb{C}\}.\]
In consequence, $\hat{\eta}(\tilde{\sigma},c^{\tilde{\sigma},G})$ is irreducible, and its equivalence class does not 
depend on the choice of $\{z_h\}_{h\in G}$. 

\section{Center construction for a C$^*$-tensor category with group action}\label{sec:Gcenter}

Let $\cD_0\subset \End_0(\cM)$ be a C$^*$-tensor category,  
$G$ a discrete  group, and $\alpha$ an action of $G$ on $\cM$ such that
$\alpha_g\circ\rho\circ \alpha_g^{-1}\in \cD_0$ for all $g\in G$, $\rho\in \cD_0$.
Namely $\alpha_g$ defines an action of $G$ on $\cD_0$. We also assume that
this action is outer, i.e., $\alpha_g\not\in \cD_0$ for $g\ne 0$.
In the following we use the notation $g[\pi]:=\alpha_g\circ \pi \circ \alpha_g^{-1}$.

Set $\cD_g:=\{\alpha_g\rho\mid \rho\in \cD_0\}$, $g\in G$, and $\cD:=\bigoplus_{g\in G}\cD_g$. 
Note that $\cD$ is also a C$^*$-tensor category.
By the outerness of $\alpha$, $\alpha_g\pi\sim \alpha_h\rho$ if and only if $g=h$, $\pi\sim \rho$.
Thus $\cD$ is  $G$-graded, and
all the results in \S \ref{sec:gradingcenter} can be applied.

Let $\cA\subset \cB$ be a Longo-Rehren inclusion for $\cD_0$.
We can extend $\alpha_g\otimes \alpha_g^{\mathrm{op}}$ to $\beta_g\in \Aut(\cB)$ by
\[
\beta_g(a\lambda_\pi):=\alpha_g\otimes \alpha_g^{\mathrm{op}}(a)\lambda_{g[\pi]}.
\]

This $\beta_g$ is nothing but the action considered in \S \ref{sec:gradingcenter}.
Moreover, $\beta_g$ preserves $\cA$ in this case, and the situation will be simplified.

\begin{df}\label{df:g-halfbraid}
Let  $\sigma\in \cD_0$,  $g\in G$.
A family of unitaries $\{\mathcal{E}^g_\sigma(\pi)\}_{\pi\in \cD_0}\subset \cM$
is called a $g$-half braiding of $\sigma$ if $\mathcal{E}_\sigma^g(\pi)\in (\sigma\pi,g[\pi]\sigma)$ and
it satisfies the condition 
$\mathrm{(BF1)}$ of Theorem \ref{thm:Gbraid1} for $\pi$.
In the following, we simply call a pair $(\sigma,\mathcal{E}^g_\sigma)$
a $g$-half braiding. 
\end{df}

Let $\{\mathcal{E}^g_\sigma(\pi)\}_{\pi\in \cD_0}$ be a $g$-half braiding for $\sigma\in \cD_0$, 
and set $\alpha_g^{-1}\cE_\sigma^g(\pi):=\alpha_g^{-1}(\cE_\sigma^g(\pi))$. 
Then
it is easy to see that $\left(\alpha_g^{-1}\sigma, \alpha_g^{-1}\mathcal{E}_\sigma^g(\cdot )\right)$
is a usual half braiding.
\begin{df}\label{df:Gcenter}
$(1)$ We define
an $\eta_G$-extension
$\eta_G(\sigma,\mathcal{E}_\sigma^g)\in \End_0(\cB)$ by
$\eta_G(\sigma,\mathcal{E}_\sigma^g):=\beta_g\eta(\alpha_g^{-1}\sigma, \alpha_g^{-1}\mathcal{E}^g_\sigma)$.
Namely, $\eta_G$ is given by
\begin{align*}
\eta_G(\sigma,\mathcal{E}_\sigma^g)(a)&=\left(\sigma\otimes \alpha_g^{\mathrm{op}}\right)(a), \,\, a\in \cA, \\
\eta_G(\sigma,\mathcal{E}_\sigma^g)(\lambda_\pi)&=\left(\mathcal{E}^g_\sigma(\pi)^*\otimes 1\right)
\lambda_{g[\pi]}, \,\, \pi\in \mathrm{Irr}(\cD_0).
\end{align*}
Similarly, we can define the opposite $\eta_G$-extension by
$\eta_G^{\mathrm{op}}(\sigma,\mathcal{E}_\sigma^g):=
\beta_g\eta^{\mathrm{op}}(\alpha_g^{-1}\sigma,\alpha_g^{-1} \mathcal{E}_\sigma^g)$.
$(2)$ The $G$-center of $(\cD_0, \alpha)$ is defined by $Z^G(\cD_0):=\bigoplus_{g\in G}Z^G(\cD_0)_g$, where   
\[
 Z^G(\cD_0)_g
:=\{\theta\in \End_0(\cB)\mid \theta\sim \eta_G(\sigma, \mathcal{E}_\sigma^g)
\mbox{ for some $g$-half braiding }(\sigma,\mathcal{E}_\sigma^g) \}.
\]
We can identify $Z_{\cD_0}(\cD)$ and $ Z^G(\cD_0)$ via identification  $\eta_G(\sigma,\cE_\sigma^g)$ with 
$ \eta(\alpha_g^{-1}\sigma,\alpha^{-1}_g\cE_\sigma^g )$.
\end{df}

We can characterize $Z^G(\cD_0)$ as follows.
\begin{align*}
Z^G(\cD_0)_g
&=\{\theta\in \End_0(\cB)\mid \theta\iota\sim \iota(\sigma \otimes \alpha_g^{\mathrm{op}})
\mbox{ for some }\sigma\in \cD_0 \} \\
&=\{\theta\in \End_0(\cB)\mid \theta \prec {\iota}(\sigma \otimes \alpha_g^{\mathrm{op}})\bar{\iota}
\mbox{ for some }\sigma\in \cD_0
\},\,\,\, (\mbox{when }|\Irr(\cD_0)|<\infty) 
\end{align*}

As in \S \ref{sec:gradingcenter}, we can see 
\[
\eta_G(\sigma,\mathcal{E}_\sigma^g)(\lambda_\pi)=\left(\mathcal{E}^g_\sigma(\pi)^*\otimes 1\right)
\lambda_{g[\pi]}
\]
holds for any (reducible) $\pi\in \cD_0$.

\bigskip

\noindent
\textbf{Remark.}
Two extensions
$\eta_G(\sigma,\mathcal{E}_\sigma^g)$ and
$\eta_G^{\mathrm{op}}(\rho,\mathcal{E}_\rho^h)$ \textit{do not commute}.
See the remark after Theorem \ref{thm:Gact}.

\bigskip

By Theorem \ref{lem:fixbraid0} (3), $\beta$ acts on $Z^G(\cD_0)$. In this case,
we can describe this action explicitly.
\begin{thm}\label{thm:Gact}
Let $(\sigma, \mathcal{E}_\sigma^g)$ be a $g$-half braiding, and
define $k[\mathcal{E}_\sigma^g](\pi):=\alpha_k\left(\mathcal{E}_\sigma^g({k^{-1}}[\pi])\right)$.
Then $(k[\sigma], k[\mathcal{E}_\sigma^g])$ is a $kgk^{-1}$-half braiding, and
$\beta_k\hateta(\sigma, \mathcal{E}_\sigma^g)\beta_k^{-1}=\hateta(k[\sigma],k[\mathcal{E}_\sigma^g])$.
\end{thm}
\textbf{Proof.}
It is straightforward to see $\{k[\mathcal{E}_\sigma^g](\pi)\}_{\pi\in \cD_0}$
is a $kgk^{-1}$-half braiding for $k[\sigma]$.
To show
$\beta_k\hateta(\sigma, \mathcal{E}_\sigma^g)\beta_k^{-1}=\hateta(k[\sigma],k[\mathcal{E}_\sigma^g])$,
we only have to verify this relation for generators of $\cB$, i.e.,
\begin{align*}
\beta_k\hateta(\sigma, \mathcal{E}_\sigma^g)\beta_k^{-1}(a)&=\hateta(k[\sigma],k[\mathcal{E}_\sigma^g])(a),\,\, a\in \cA, \\ \beta_k\hateta(\sigma, \mathcal{E}_\sigma^g)\beta_k^{-1}(\lambda_\pi)&=\hateta(k[\sigma],k[\mathcal{E}_\sigma^g])(\lambda_\pi),
\,\, \pi\in \mathrm{Irr}(\cD_0).  
\end{align*}
These can be seen by simple computation.
\hfill$\Box$

\bigskip

\noindent
\textbf{Remark.}
Remarked as above,
$\eta_G(\sigma,\mathcal{E}_\sigma^g)$ and
$\eta_G^{\mathrm{op}}(\rho,\mathcal{E}_\rho^h)$ do not commute.
Note two endomorphisms
\[
\beta_g^{-1}\eta_G(\sigma,\mathcal{E}_\sigma^g)
=\eta(\alpha_g^{-1}\sigma,\alpha_g^{-1} \mathcal{E}_\sigma^g),\,\,
\beta_h^{-1}\eta_G^{\mathrm{op}}(\rho,\mathcal{E}_\rho^h)=\eta^{\mathrm{op}}(\alpha_h^{-1}\sigma,
\alpha_h^{-1} \mathcal{E}_\rho^h)
\]
commute.
So
\[
\beta_g^{-1}\eta_G(\sigma,\mathcal{E}_\sigma^g) \beta_h^{-1}\eta_G^{\mathrm{op}}(\rho,\mathcal{E}_\rho^h)
=
\beta_h^{-1}\eta_G^{\mathrm{op}}(\rho,\mathcal{E}_\rho^h)\beta_g^{-1}\eta_G(\sigma,\mathcal{E}_\sigma^g).
\]
Hence we have the following commutation relation by applying Theorem \ref{thm:Gact};
\[
\beta_{hg}^{-1}\eta_G(h[\sigma],h[\mathcal{E}_\sigma^g]) \eta_G^{\mathrm{op}}(\rho,\mathcal{E}_\rho^h)
=
\beta_{gh}^{-1}\eta_G^{\mathrm{op}}(g[\rho],g[\mathcal{E}_\rho^h])\eta_G(\sigma,\mathcal{E}_\sigma^g).
\]

\bigskip

As in \cite{Iz-LR}, \cite{Kw-Dcenter}, we will relate $Z^G(\cD_0)$ with a Tube algebra associated with $\cD_0$ and $\alpha$.
\begin{df}\label{df:G-tube}
We define a $G$-twisted Tube algebra $\mathrm{Tube}^G(\cD_0)$ for $\cD_0$ as follows.
As a linear space,
\[
\mathrm{Tube}^G(\cD_0)_g:=\bigoplus_{\substack{\pi,\rho, \sigma\in \mathrm{Irr}(\cD_0)
}}
\left(\sigma \pi, g[\pi]\rho\right),\,\,\,
\mathrm{Tube}^G(\cD_0):=\bigoplus_{g\in G}\mathrm{Tube}^G(\cD_0)_g.
\]
Define a product and $*$-operation as follows. (We use the notation of \cite{Iz-LR}.)
\begin{align*}
\langle\sigma \pi|X|g[\pi]\rho \rangle
\langle\sigma' \pi'|Y|h[\pi']\rho' \rangle&\!:=\delta_{g,h}\delta_{\rho,\sigma'}
\bigoplus_{\xi\in \mathrm{Irr}(\cD_0)}\sum_{e}
\langle\sigma \xi|\alpha_g(T_{\pi,\pi'}^{\xi,e*})g[\pi](Y)X \sigma(T_{\pi,\pi'}^{\xi,e})
|g[\xi]\rho' \rangle, \\
\langle\sigma \pi|X|g[\pi]\rho \rangle^*&\!:=
\langle\rho \bar{\pi} | g[\bar{\pi}]\left(\sigma(\bar{R}_\pi)X^*\right)\alpha_g(R_\pi)|g[\bar{\pi}]
\sigma \rangle.
\end{align*}
Here $\{T_{\pi,\pi'}^{\xi,e}\}_e\subset (\xi,\pi\pi')$ is an orthonormal basis.
\end{df}

\begin{thm}\label{thm:GTubecenter}
Suppose $[\cB:\cA]<\infty$. 
There exists 1 to 1 correspondence between $\mathrm{Irr}(Z^G(\cD_0)_g)$ and the set of
minimal central projections of $\mathrm{Tube}^G(\cD_0)_g$.
\end{thm}
\textbf{Proof.}
It is easy to see the canonical isomorphism
\[
(\alpha_g^{-1}\rho\pi, \pi \alpha_g^{-1}\rho)\ni T\rightarrow \alpha_g(T)\in
(\rho\pi, g[\pi] \rho)
\]
of Hilbert spaces gives an $*$-isomorphism of $\mathrm{Tube}(\cD_0,\cD)$ and
$\mathrm{Tube}^G(\cD_0)$.
Together with Lemma \ref{lem:diffgrade}, we get the conclusion.
\hfill$\Box$

\bigskip

Proposition \ref{prop:Vunitary}
implies the following.
\begin{thm}\label{thm:hatetaop}
Let $(\sigma,\mathcal{E}_\sigma^g)$ be a $g$-half braiding, and 
define
\[
\mathcal{E}_{\bar{\sigma}}^{g^{-1}}(\pi):=R_\sigma^*\bar{\sigma}(\mathcal{E}_\sigma(g^{-1}(\pi))^*)\bar{\sigma}\pi(\bar{R}_\sigma).
\]
Then we have
$\overline{\eta_G(\sigma, \mathcal{E}_\sigma^{g})}\sim 
\eta_G(\bar{\sigma}, \mathcal{E}_{\bar{\sigma}}^{g^{-1}})
\sim\beta_g^{-1}\hateta(\sigma, \mathcal{E}_\sigma^g)\sim\hateta(\sigma, \mathcal{E}_\sigma^g)\beta_g^{-1}$.
\end{thm}
\textbf{Proof.}
We can see $\overline{\eta_G(\sigma, \mathcal{E}_\sigma^g)}\sim 
\eta_G(\bar{\sigma}, \mathcal{E}_{\bar{\sigma}}^{g^{-1}})$ in a similar way as Proposition \ref{prop:eta-inter} (2).
By 
Proposition \ref{prop:Vunitary}
\[
\overline{ \eta_G^{\mathrm{op}}(\sigma, \mathcal{E}_\sigma^g)}=
\overline{\beta_g\eta^{\mathrm{op}}(\alpha_g^{-1}\sigma, \alpha_g^{-1}\mathcal{E}_\sigma^g)}\sim
\eta(\alpha_g^{-1}\sigma, \alpha_g^{-1}\mathcal{E}_\sigma^g)
=\beta_g^{-1} \eta_G(\sigma, \mathcal{E}_\sigma^g).
\]
Note that $\alpha_g^{-1}\sigma$ is in a $g^{-1}$-component of $\cD$.
\hfill$\Box$

\bigskip

At the end of this section, we will show that $Z^G(\cD_0)$ has a $G$-braiding.
The following description of intertwiner spaces 
can be proved in the same way as in \cite[Theorem 4.6]{Iz-LR}, \cite[Theorem 2.4]{BEK-LR}.
\begin{thm}\label{thm:GetaIntertwiner}
Let $\eta_G(\sigma, \mathcal{E}_\sigma^g), \eta_G(\rho, \mathcal{E}_{\rho}^{h})\in Z^G(\cD_0)$.
Then $\left(\eta_G(\sigma, \mathcal{E}_\sigma^g), \eta_G(\rho, \mathcal{E}_\rho^{h})\right)=0$
for $g\ne h$. When $g=h$,
\[
\left(\eta_G(\sigma, \mathcal{E}_\sigma^g), \eta_G(\rho, \mathcal{E}_\rho^{g})\right)
=\{T\otimes 1\mid T \in (\sigma,\rho),
\mathcal{E}_{\rho}^g(\pi)T=g[\pi](T)\mathcal{E}_\sigma^g(\pi),\pi\in \cD_0\}.
\]
\end{thm}

\bigskip

In a similar way as in \cite[Theorem 4.6 (iii)]{Iz-LR}, we can show the existence of a $G$-braiding.
In contrast to Theorem \ref{thm:Gbraid1}, it is easy to see the properties of a $G$-braiding.
\begin{thm}\label{thm:Gbraid2}
Define $\mathcal{E}_G(\tilde{\sigma},\tilde{\rho}):=\mathcal{E}_\sigma^g(\rho)\otimes 1$ 
for
$\tilde{\sigma}=\eta_G(\sigma, \mathcal{E}_\sigma^g)\in Z^G(\cD_0)_g$ and
$\tilde{\rho}=\eta_G(\rho, \mathcal{E}_\rho^h)\in Z^G(\cD_0)_h$.
Then $\mathcal{E}_G$ gives a $G$-braiding on $Z^G(\cD_0)$. 
\end{thm}
\textbf{Proof.}
First note that
endomorphisms
$\tilde{\sigma}\tilde{\rho}$ and $g[\tilde{\rho}]\tilde{\sigma}$
are given by
$\tilde{\sigma}\tilde{\rho}=\eta_G(\sigma\rho, \mathcal{E}_{\sigma\rho}^{gh})$ and
$g(\tilde{\rho})\tilde{\sigma}=\eta_G\left(g[\rho]\sigma, \mathcal{E}_{g[\rho]\sigma}^{gh}\right)$,
respectively.
Here
\[
\mathcal{E}_{\sigma\rho}^{gh}(\pi)=\mathcal{E}_\sigma^g(h(\pi))\sigma(\mathcal{E}^h_\rho(\pi)),\,
\mathcal{E}_{g[\rho]\sigma}^{gh}(\pi)=g[\mathcal{E}_\rho^h](g[\pi])
g[\rho](\mathcal{E}^g_\sigma(\pi))
=\alpha_g(\mathcal{E}_\rho^h(\pi))g[\rho](\mathcal{E}_\sigma^g(\pi)).
\]
To show
$\mathcal{E}_G(\tilde{\sigma},\tilde{\rho})\in (\tilde{\sigma}\tilde{\rho}, g[\tilde{\rho}]\tilde{\sigma})$,
we only have to show
\[
\alpha_g(\mathcal{E}_\rho^h(\pi))g[\rho](\mathcal{E}_\sigma^g(\pi))\mathcal{E}_\sigma(\rho)=
gh[\pi]\left(\mathcal{E}_\sigma^g(\rho)\right)\mathcal{E}_\sigma^g(h[\pi])\sigma(\mathcal{E}^h_\rho(\pi))
\]
by Theorem \ref{thm:GetaIntertwiner}.
By the half braiding property  of $\mathcal{E}_\sigma^g$, the above equation is equivalent to
the following.
\[
\alpha_g(\mathcal{E}_\rho^h(\pi))\mathcal{E}^g_\sigma(\rho\pi)=
\mathcal{E}_\sigma^g(h[\pi]\rho)\sigma(\mathcal{E}^h_\rho(\pi)).
\]
Since $h[\pi]\rho= \Ad \left(\mathcal{E}_\rho^h(\pi)\right)\circ \rho\pi$,
we have
$\mathcal{E}_\sigma^g(h[\pi]\rho)=\alpha_g\left(\mathcal{E}_\rho^h(\pi)\right)
\mathcal{E}_\sigma^g(\rho\pi)\sigma\left(\mathcal{E}_\rho^h(\pi)\right)^*$.
Thus the above equation actually holds.

The condition (BF1) for $\mathcal{E}_G(\tilde{\sigma},\tilde{\rho})$
follows from the half braiding condition of $\mathcal{E}_\sigma^g$.
The  condition (BF2) for $\mathcal{E}_G(\tilde{\sigma},\tilde{\rho})$
follows from the description of intertwiner spaces
between $\eta_G$-extension in Theorem \ref{thm:GetaIntertwiner}.
As $k[\tilde{\sigma}]=\hateta(k[\sigma],k[\mathcal{E}_\sigma^g])$,
\[
\mathcal{E}_G\left(k[\tilde{\sigma}], k[\tilde{\rho}]\right)=k[\mathcal{E}_\sigma^g](k[\rho])\otimes 1=
\alpha_k(\mathcal{E}_\sigma^g(\rho))\otimes 1=\beta_k(\mathcal{E}_G(\tilde{\sigma},\tilde{\rho})).
\]
So we obtain the condition (BF3).
\hfill$\Box$

\bigskip

\noindent
\textbf{Remark 1.} In \S \ref{sec:gradingcenter}, we construct a $G$-braiding $\mathcal{E}$ on $Z_{\cD_0}(\cD_0)$.
So we can apply this construction to $Z^G(\cD_0)=Z_{\cD_0}(\cD)$.
We can see 
\[
\mathcal{E}_G(\tilde{\sigma},\tilde{\rho})=
\beta_{gh}\left(\mathcal{E}^-\left(h^{-1}[\widetilde{\alpha_g^{-1}\sigma}], \widetilde{\alpha_h^{-1}\rho}\right)\right)=
\beta_{gh}\left(\mathcal{E}^+(\widetilde{\alpha_h^{-1}\rho},\widetilde{\alpha_g^{-1}\sigma})^*\right).
\]

\medskip

\noindent
\textbf{Remark 2.}
By combining methods in \S \ref{sec:gradingcenter} and \S \ref{sec:Gcenter},
it may be possible to  present Oikawa's center construction \cite{Oikawa-center} by using Longo-Rehren inclusion.

\appendix

\section{Realization of a Tube algebra in a Longo-Rehren inclusion}\label{sec:TubeReal}

In this appendix, we explain the detail of M{\"u}ger's remark in \cite[Remark 5.1]{Mu-subII}
by realizing a Tube algebra in a Longo-Rehren inclusion explicitly.

We use the notation in \S \ref{sec:gradingcenter}.
Fix a family $\{w_\xi \}_{\xi \in \cC_0}\subset \cA$ of orthogonal isometries with support 1, 
and define a unitary $W\in B(\fH_0, L^2(\cA))$ by 
$W v:=\sum_{\xi\in \Irr(\cC_0)}w_\xi v(\xi)$.
It is clear that $W$ commutes with the right action of $\cA$, and 
hence the conjugate $\bar{\iota}\in \End(\cB_0, \cA)$ of the inclusion map 
$\iota:\cA\rightarrow \cB_0$ is given by $W \gamma_L(a)W^*$ for $a\in \cB_0$.
We can see 
\[
\bar{\iota}(a)=\sum_{\xi\in \cC_0} w_\xi (\xi\otimes \xi^{\mathrm{op}})(a)w_\xi^*, \,\,a\in \cA, \,\,
\bar{\iota}(\lambda_\eta)=\sum_{\xi,\zeta\in \cC_0}\sqrt{\dfrac{d(\xi)}{d(\zeta)}}w_\xi \tilde{T}_{\xi,\eta}^\zeta w_\zeta^*.
\]

\noindent
\textbf{Remark.} 
We can show that  
$\sum\limits_{\xi\in \cC_0} (s_\xi^0\otimes \bar{R}_\xi^*)\lambda_{{\xi}} 
(\mu\otimes \id)(w_{\bar{\xi}}^{*}) $
 is a unitary in $(\iota (\mu\otimes \id)\bar{\iota}, \tilde{\theta}_\mu^{(0)})$
by using above equations. 

\bigskip

Let $\kappa_\xi^0:=\sum\limits_{\eta,\zeta\in \cC_0}\sqrt{\dfrac{d(\eta)}{d(\zeta)}}\xi\otimes \xi^\mathrm{op}(w_\eta)\tilde{T}_{\xi,\eta}^\zeta w_\zeta^*$.
By using above facts,  we can verify $\kappa_\xi^0 \in (\bar{\iota},(\xi\otimes \xi^{\mathrm{op}})\bar{\iota})$.

For $T\in (\pi\xi,\xi\rho)$, $\pi,\rho\in \cC$, $\xi\in \cC_0$, 
set 
\[
\theta^{(l)}(T):=\lambda_\xi^*(T\otimes 1)(\pi\otimes \id)(\kappa_\xi^0), \,\,\,
\theta^{(r)}(T):=(\id\otimes \pi^{\mathrm{op}})(\kappa_\xi^{0*})
(1\otimes j(T^*))\lambda_\xi.\]
It is easy to see 
\[
 \theta^{(l)}(T)\in (\iota(\pi\otimes \id)\bar{\iota},\iota(\rho\otimes \id)\bar{\iota}),\,\,\,
 \theta^{(r)}(T)\in (\iota(\id\otimes \rho^{\mathrm{op}})\bar{\iota},
\iota(\id\otimes \pi^{\mathrm{op}})\bar{\iota}).
\]

Let us  denote the product and the 
$*$-operation of $\mathrm{Tube}(\cC_0,\cC)$ by $T\bullet S$  and $T^{\star}$ for 
$T\in (\pi\xi,\xi\rho)$, $S\in (\rho\eta,\eta \sigma)$.
(We adopt the definition of Tube algebras given in \cite{Iz-LR}.)
Then  
\begin{align*}
 \theta^{(l)}(S)\theta^{(l)}(T)&=
\lambda_\eta^*(S\otimes 1)(\rho\otimes \id)(\kappa_\eta^0)\lambda_\xi^*(T\otimes 1)(\pi\otimes \id)(\kappa^0_\xi)\\
&=
\lambda_\eta^*\lambda_\xi^*(\xi(S)\otimes 1)\xi\rho\otimes \xi^{\mathrm{op}}(\kappa^0_\eta)
(T\otimes 1)(\pi\otimes \id)(\kappa^0_\xi)\\
&=
\lambda_\eta^*\lambda_\xi^*(\xi(S)T\otimes 1)\pi\xi\otimes \xi^{\mathrm{op}}(\kappa^0_\eta)
(\pi\otimes \id)(\kappa^0_\xi)\\
&=\sum_{\zeta,\zeta'}
\lambda_\zeta^*T_{\xi,\eta}^{\zeta*}(\xi(S)T\otimes 1)\pi(T_{\xi,\eta}^{\zeta'})(\pi\otimes \id)(\kappa^0_\pi) \\
&=\sum_{\zeta,e}
\lambda_\zeta^*(T_{\xi,\eta}^{\zeta,e*}\xi(S)T\pi(T_{\xi,\eta}^{\zeta,e})\otimes 1)(\pi\otimes \id)(\kappa^0_\pi) \\
&=\sum_{\zeta,e}
\theta^{(l)}(T_{\xi,\eta}^{\zeta,e*}\xi(S)T\pi(T_{\xi,\eta}^{\zeta,e})) 
=\theta^{(l)}(T\bullet S),
\end{align*}
and 
\begin{align*}
 \theta^{(l)}(T)^*&=
(\pi\otimes \id)(\kappa^0_\xi)^*(T^*\otimes 1)\lambda_\xi  
=(\pi\otimes \id)(\bar{R}_{\xi\otimes \xi^{\mathrm{op}}}^*)
(\pi\xi\otimes \xi^{\mathrm{op}})(\kappa^0_{\bar{\xi}})(T^*\otimes 1)\lambda_\xi \\
&=
(\pi\otimes\id )(\bar{R}_{\xi\otimes \xi^{\mathrm{op}}}^*)(T^*\otimes 1)(\xi\rho\otimes 
\xi^{\mathrm{op}})(\kappa^0_{\bar{\xi}})
\lambda_\xi \\
&=
(\pi\otimes \id)(\bar{R}_{\xi\otimes \xi^{\mathrm{op}}}^*)(T^*\otimes 1)
\lambda_\xi 
\rho(\kappa^0_{\bar{\xi}}) \\
&=
(\pi\otimes \id)(\bar{R}_{\xi\otimes \xi^{\mathrm{op}}}^*)(T^*\otimes 1)
\lambda_{\bar{\xi}}^*{R}_{\xi\otimes \xi^{\mathrm{op}}}
\rho(\kappa^0_{\bar{\xi}})  \\
&=
\lambda_{\bar{\xi}}^*(
\bar{\xi}\pi\otimes \bar{\xi}^{\mathrm{op}})(\bar{R}_{\xi\otimes \xi^{\mathrm{op}}}^*)(\bar{\xi}(T)^*\otimes 1)
{R}_{\xi\otimes \xi^{\mathrm{op}}}
\rho(\kappa^0_{\bar{\xi}}) \\
&=
\lambda_{\bar{\xi}}^*
(\bar{\xi}\pi(\bar{R}_\xi^*)
\bar{\xi}(T)^*{R}_{\xi} \otimes 1)
\rho(\kappa^0_{\bar{\xi}}) 
=\theta^{(l)}(T^{\star})
\end{align*}
holds.

Fix a family of orthogonal isometries$\{v_\pi\}_{\pi\in \cC}\subset \cB_0$ with
$\sum_{\pi\in \cC}v_\pi v_\pi^*=1$.  
Let $A^{(l)}$ and  $A^{(r)}$  
 be  linear span of $\{v_\rho \theta^{(l)}(T) v_\pi^*\mid T\in (\pi\xi,\xi\rho)\}$
and 
$\{v_\pi \theta^{(r)}(T) v_\rho^*\mid T\in (\pi\xi,\xi\rho)\}$, respectively.
By the above computation, 
$\Theta^{(l)}: T\in (\pi\xi,\xi\rho)\rightarrow v_\rho\theta^{(l)}(T)v_\pi^*\in A^{(l)}$ 
is an anti-isomorphism 
of $\mathrm{Tube}(\cC_0,\cC)$ onto $A^{(l)}$. 
In a similar way, we can verify that
$\Theta^{(r)}: T\in (\pi\xi,\xi\rho)\rightarrow v_\pi\theta^{(r)}(T)v_\rho^*\in A^{(r)}$
is an isomorphism of  
$\mathrm{Tube}(\cC_0, \cC)$ onto $A^{(r)}$. 

Let $a\in (\iota(\pi\otimes \id)\bar{\iota},\iota(\rho\otimes \id)\bar{\iota})$,
and expand $a=\sum_{\xi \in \cC_0 }\lambda_\xi^*a_\xi$. Then
$a_\xi\in ((\pi\otimes \id)\bar{\iota}, (\xi\rho\otimes \xi^{\mathrm{op}} )\bar{\iota})$.
Hence we have two conditions
\[
a_\xi\in ((\pi\otimes \id)\bar{\iota}\iota, (\xi\rho\otimes \xi^{\mathrm{op}}) \bar{\iota}\iota),\,\,
a_\xi (\pi\otimes \id)\bar{\iota}(\lambda_\eta)=\xi\rho\otimes \xi^{\mathrm{op}} \bar{\iota}(\lambda_\eta)a_\xi,\eta\in \cC_0.
\]
From the former condition, we obtain
$\rho\xi\otimes \xi^{\mathrm{op}}(w_0^*)
 a_\xi \pi\otimes \id(w_\eta)
\in (\pi \eta\otimes \eta^{\mathrm{op}}, \xi\rho\otimes \xi^{\mathrm{op}})$.
Then 
$ \rho\xi\otimes \xi^{\mathrm{op}}(w_0^*)
 a_\xi \pi\otimes \id(w_\eta)=0$ 
for $\xi\ne \eta$, 
and 
$ \rho\xi\otimes \xi^{\mathrm{op}}(w_0^*) a_\xi \pi\otimes \id(w_\xi)=T_{\xi}\otimes 1$ for some 
$T_{\xi}\in (\pi\xi,\xi\rho)$,
and 
\[
 \rho\xi\otimes \xi^{\mathrm{op}}(w_0^*)
 a_\xi=
\left(\rho\xi\otimes \xi^{\mathrm{op}}(w_0^*)
 a_\xi \right)\left(\sum_{\eta\in \cC_0}\pi\otimes \id(w_\eta w_\eta^*)\right)=(T_\xi\otimes 1)\pi\otimes \id(w_\xi^*).
\]

Multiplying $ \rho\xi\otimes \xi^{\mathrm{op}}(w_0^*)$ from the left side of 
the latter condition, we have
\begin{align*}
& \rho\xi\otimes \xi^{\mathrm{op}}(w_0^*)
 a_\xi \pi\otimes \id
(\bar{\iota}(\lambda_\eta))=
\xi\rho\otimes \xi^{\mathrm{op}} (w_0^*\bar{\iota}(\lambda_\eta))a_\xi  \\
&\Rightarrow 
 \xi\rho\otimes \xi^{\mathrm{op}} (w_\eta^*)a_\xi =
(T_\xi \otimes 1) \pi\otimes \id
(w_\xi ^*\bar{\iota}(\lambda_\eta))=
\sum_{\zeta\in \cC_0}(T_\xi \otimes 1) \pi\otimes \id
(\tilde{T}_{\xi,\eta}^{\zeta}w_\zeta^*).
\end{align*}
Hence we get
\begin{align*}
 a_\xi &= \sum_{\eta\in \cC_0}\xi\rho\otimes \xi^{\mathrm{op}} (w_\eta w_\eta^*)a_\xi  
=
\sum_{\eta,\zeta\in \cC_0}\xi\rho\otimes \xi^{\mathrm{op}} (w_\eta )
(T_\xi \otimes 1) \pi\otimes \id
(\tilde{T}_{\xi,\eta}^{\zeta}w_\zeta^*) \\
&=
\sum_{\eta,\zeta\in \cC_0}(T_\xi \otimes 1) 
\pi\xi \otimes \xi^{\mathrm{op}} (w_\eta )
\pi\otimes \id
(\tilde{T}_{\xi,\eta}^{\zeta}w_\zeta^*) 
=(T_\xi \otimes  1) \pi\otimes \id (\kappa_\xi^0).
\end{align*}
So we conclude 
$a\in (\iota(\pi\otimes \id)\bar{\iota},\iota(\rho\otimes \id)\bar{\iota})$ is 
of the form $a=\sum_{\xi\in \cC_0} \theta^{(l)}(T_\xi)$, $T_\xi \in (\pi\xi,\xi\rho)$,
and $A^{(l)}$ is dense in $\sum_{\pi,\rho\in \cC}v_\rho 
(\iota(\pi\otimes \id)\bar{\iota},\iota(\rho\otimes \id)\bar{\iota})v_\pi^*$.
(If $|\Irr(\cC)|<\infty$, we have
$A^{(l)}=\sum_{\pi,\rho\in \cC}v_\rho 
(\iota(\pi\otimes \id)\bar{\iota},\iota(\rho\otimes \id)\bar{\iota})v_\pi^*$.)

With a little bit of effort, we can present a result of \cite{De-Kac} in terms 
of Longo-Rehren inclusions by applying results in this appendix,
c.f. \cite{Nesh-Yama-remtube}.


\ifx\undefined\bysame
\newcommand{\bysame}{\leavevmode\hbox to3em{\hrulefill}\,}
\fi

\end{document}